\documentclass[11pt]{article}
\usepackage{amsmath}
\usepackage{amsfonts}
\usepackage{color}
\usepackage{hyperref}
\usepackage{dpfloat}

\usepackage[pdftex]{graphicx}

\newtheorem{theo}{Theorem}[section]
\newtheorem{lem}[theo]{Lemma}
\newtheorem{cor}[theo]{Corollary}

\newcommand{\mysection}[1]{\section{#1} \setcounter{equation}{0}}
\newcommand{\proof}{{\sc Proof.} \quad}
\newcommand{\proofc}{{\sc Proof} \ }
\newcommand{\be}{\begin{equation} \label}
\newcommand{\ee}{\end{equation}}
\newcommand{\bea}{\begin{eqnarray}\label}
\newcommand{\eea}{\end{eqnarray}}
\newcommand{\bas}{\begin{eqnarray*}}
\newcommand{\eas}{\end{eqnarray*}}
\newcommand{\bit}{\begin{itemize}}
\newcommand{\eit}{\end{itemize}}
\newcommand{\qed}{\hfill$\Box$ \vskip.2cm}
\newcommand{\nn}{\nonumber}
\newcommand{\R}{\mathbb{R}}
\newcommand{\N}{\mathbb{N}}
\newcommand{\pO}{\partial\Omega}
\newcommand{\bom}{\overline{\Omega}}
\newcommand{\eps}{\varepsilon}

\newcommand{\hra}{\hookrightarrow}
\newcommand{\io}{\int_\Omega}

\newcommand{\abs}{\\[5pt]}
\newcommand{\tm}{T_{max}}
\newcommand{\Du}{D_u}
\newcommand{\Dh}{D_h}
\newcommand{\Cf}{C_f}
\newcommand{\Cg}{C_g}
\newcommand{\cphi}{c_\phi}
\newcommand{\Cphi}{C_\phi}
\newcommand{\Cpsi}{C_\psi}
\newcommand{\CPhi}{C_\Phi}
\newcommand{\hatt}{\widehat{T}}
\newcommand{\F}{{\mathcal{F}}}
\newcommand{\D}{{\mathcal{D}}}
\begin{document}
\enlargethispage{10mm}
\title{Does indirectness of signal production reduce the explosion-supporting potential in chemotaxis-haptotaxis systems?\\
Global classical solvability in a class of models for cancer invasion (and more)}
\author{
Christina Surulescu\footnote{surulescu@mathematik.uni-kl.de}\\
{\small Technische Universit\"{a}t Kaiserslautern, Felix-Klein-Zentrum f\"{u}r Mathematik,} \\
{\small 67663 Kaiserslautern, Germany}
\and
Michael Winkler\footnote{michael.winkler@math.uni-paderborn.de}\\
{\small Institut f\"ur Mathematik, Universit\"at Paderborn,}\\
{\small 33098 Paderborn, Germany} 
}

\date{}
\maketitle
\begin{abstract}
\noindent 
We propose and study a class of parabolic-ODE models involving chemotaxis and haptotaxis of a species following signals indirectly produced by another, non-motile one. 
The setting is motivated by cancer invasion mediated by 
interactions with the tumor microenvironment, but has much wider applicability, being able to comprise descriptions of biologically quite different problems.
As a main mathematical feature consituting a core difference to both classical Keller-Segel chemotaxis systems
and Chaplain-Lolas type chemotaxis-haptotaxis systems, the considered model accounts for certain types of indirect
signal production mechanisms.\abs
The main results assert unique global classical solvability under suitably mild assumptions on the system parameter
functions in associated spatially two-dimensional initial-boundary value problems.
In particular, this rigorously confirms that at least in two-dimensional settings, the con\-si\-dered indirectness in
signal production induces a significant blow-up suppressing tendency also in taxis systems substantially more general than
some particular examples for which corresponding effects have recently been observed.\abs
\end{abstract}
%
%
%
%
%
%
%
%
\section{Introduction}\label{intro}
We study here a general class of parabolic-parabolic-ODE-ODE systems (see \eqref{0} below) containing the following 
model of cancer invasion with chemotaxis and haptotaxis:
\be{02}
	\left\{ \begin{array}{l}
	u_t=\Delta u - \chi \nabla \cdot (u\nabla h) - \xi \nabla \cdot (u\nabla v) + \mu u(1-u-v-w), \\[1mm]
	h_t=\Delta h - h + \alpha w, \\[1mm]
	v_t=-hv + \eta v(1-u-v)+\beta \frac{w}{1+w}, \\[1mm]
	w_t=u,
	\end{array} \right.
\ee
supplemented with adequate initial conditions and no-flux boundary conditions. 
The model variables are: $u$: density of tumor cells, $h$: concentration of matrix metalloproteinases (MMP), $v$: 
density of tissue fibers (extracellular matrix, ECM), $w$: density of cancer associated fibroblasts (CAFs). 
The chemotactic bias of the cells is in the direction of the MMP gradient, while haptotaxis means as usual following the
gradient of tissue density. In most of the previous chemotaxis and chemotaxis-haptotaxis models the chemoattractant is directly
produced by the population performing diffusion and taxis, thereby involving (\cite{SSW}) or not 
(\cite{pang_wang_JDE}, \cite{tao_wang_2009}, \cite{taowin_JDE2014}) the ECM in this production. In our present model the
chemoattractant is generated in an indirect way, by the CAFs,  which are activated by the tumor cells. In \cite{chaplain-lolas}
was introduced a complex model for the evolution of a population of tumor cells interacting with two chemoattractants and also
performing haptotaxis. One of the chemoattractants therein is directly produced by the cells, while the other's production is
mediated by another substance, the latter being in turn produced by the cells under the influence of the first chemoattractant.
All model variables (except tissue) are diffusing in a linear way; in particular, the producers of all chemoattractants are
diffusive. This is not the case in our model \eqref{02}; moreover, both signals are completely ($h$) or partially ($v$)
obtained from the non-diffusing, indirect producer $w$.\abs
Model \eqref{02} is motivated by the problem of investigating CAF-mediated cancer invasion into the surrounding tissue. 
CAFs are major components of the neoplasm microenvironment. They secrete a variety of extracellular matrix components and are
involved in the formation of the desmoplastic stroma  characterizing many advanced carcinomas (\cite{hanahan-w}). 
For a long time the general belief was that tumor development, invasion, and metastasis occur 
as a result of cancer progression. Recent studies revealed, however, that CAFs contribute instead of tumor cells to these
processes, via expression of various growth factors, cytokines, chemokines, and degradation of ECM (\cite{cirri}, 
\cite{hanahan-w}, \cite{quail}, \cite{spano}), but also by restructuring the latter to facilitate migration (\cite{jolly}).
Fibronectin (Fn) assembled by CAFs mediates cell association and directional migration. Compared with normal fibroblasts, CAFs
produce an Fn-rich ECM with anisotropic fiber orientation, along which the tumor cells preferentially migrate (\cite{erdogan}).
 The origin of CAFs is not completely elucidated; we refer e.g. to \cite{cirri-chiarugi-2,ishii} for a couple of reviews. 
There  is evidence that they can arise among others from carcinoma cells through 
epithelial-mesenchymal transition (EMT) (\cite{kalluri}), thus allowing the cancer cells to adopt a mesenchymal phenotype associated to enhanced migratory capacity and invasiveness (\cite{duda}, \cite{neri}). 
It has also been shown (cf.~e.g.~\cite{mitra}) that cancer cells can reprogram resident tissue fibroblasts to become CAFs 
through the actions of miRNAs. MMPs are primarily
derived from CAFs in various types of tumor (\cite{kalluri}, \cite{lecomte}, \cite{tao-caf}). 
In particular, it has been shown e.g., that MMP-9, an endoproteinase involved in ECM degradation and implicated as a
prerequisite of metastasis, has very limited or no expression in various cancer cell lines. Instead, MMP-9 is well-known to be
secreted from cancer stromal fibroblasts and endothelial cells (\cite{stuelten}). 
For further information about CAFs and their MMP production we refer, e.g., to the review \cite{lecomte}.
The main features of CAF-mediated tumor invasion mentioned above are captured in our model \eqref{02} below.\abs
%
%
Closely related from a mathematical viewpoint is the chemotaxis-haptotaxis model 
\be{01}
\left\{ \begin{array}{l}
         u_t=D_u\Delta u+\chi \nabla \cdot (u\nabla h)-\xi \nabla \cdot (u\nabla v)-k_1u+k_2\frac{hw}{1+hw},\\[1mm]
         h_t=D_h\Delta h+k_3w-k_4h,\\[1mm]
         v_t=-k_5(h+u+w)v+k_6v(1-v)_+^2,\\[1mm]
         w_t=k_7u-k_8hw+k_9w(1-u-v-w)_+^2,
        \end{array}\right .
        \ee
which describes the evolution of two cancer cell subpopulations either proliferating ($w$) or  migrating ($u$), with the
corresponding transitions between the two phenotypes, and in interaction with tissue ($v$) and acidity ($h$). The latter is
mainly produced by the highly glycolytic proliferating tumor cells -- hence indirect signal production. 
The proliferation/migration (also known as go-or-grow) dichotomy asserting that moving cells defer their proliferation seems to
be a relevant feature for some types of cancer (\cite{giese}, \cite{jerby},\cite{widmer}) and can lead to interesting 
mathematical problems and qualitative behavior of the corresponding heterogeneous tumor (\cite{SSU}, \cite{ZSH}). 
The migrating cells in \eqref{01} perform haptotaxis and pH-taxis: they move towards increasing ECM gradient and away from
acidic (thus hypoxic) areas. The conversion from proliferating to migrating cells depends on the concentration $h$ of protons
in the peritumoral region and infers limitations, as only a rather small part of the tumor becomes motile -usually cells
situated at the tumor margins. Moreover, the protons are buffered by the environment (e.g., uptake by vasculature), contribute
to tissue degradation, and restrict tumor proliferation. Supplementary to hypoxia the tissue can be degraded by chemical or
biological agents directly or indirectly produced or activated by the tumor cells, see e.g. \cite{madsen} and \cite{matus}. 
The ECM remodeling also involves limited, logistic-like growth. The latter is also used in \eqref{01} to describe tumor cell proliferation. We note here that taking squares of the positive parts of the proliferation terms is motivated by the 
mainly technical need of satisfying the regularity assumptions in Theorem \ref{theo22} and could perhaps be relaxed;
in fact, the difference between functions of this type just involving the positive parts and those taking their squares
is rather small. System \eqref{01} fits in the theoretical framework proposed and analyzed below.\abs
Thus the general model class to be introduced and studied here includes descriptions of biologically quite different problems, as exemplified above. It therefore provides a comprehensive mathematical structure 
for several issues related to cancer cell migration under the influence of chemo- and haptotactic effects, including 
some variants of the often addressed Chaplain-Lolas model from \cite{chaplain-lolas}. 
Other problems characterizing tactic cell migration, e.g. in wound healing and/or angiogenesis or interacting microbial 
populations in biofilm formation and persistence, can potentially be cast 
into this mathematical framework. Models describing the dynamics of a species performing chemotaxis towards the gradient of a 
signal produced by another, non-tactic species 
(see e.g.~\cite{hutao} and \cite{taowin_JEMS}) form a subclass of this setting.\\[3mm]
{\bf Approaching a mathematical core feature: indirectness of chemoattractant production.} \quad
From a purely mathematical perspective, 
a common feature distinguishing both (\ref{01}) and (\ref{02}) from classical chemotaxis or
chemotaxis-haptotaxis systems of Keller-Segel or Chaplain-Lolas type consists in the circumstance that the respective
mechanisms of chemoattractant production are indirect in the sense that the corresponding signal is not produced
directly by individuals of the cell population, but rather through a third agent.
Possible implications on the system dynamics, however, have apparently been detected only in some quite particular examples
of chemotaxis-only models: 
Indeed, the only findings in this direction we are aware of concentrate on associated derivates of the classical
Keller-Segel system, for which substantial blow-up preventing effects of such indirect signal production mechanisms
have recently been revealed. More precisely, in sharp contrast to what is known for 
standard Keller-Segel systems (\cite{horstmann_DMV}, \cite{BBTW}), 
regardless of the size of the initial data some 
corresponding spatially two-dimensional initial-boundary value problems always possess globally
defined classical solutions (\cite{taowin_JEMS}, \cite{laurencot}). \abs 
Complementary to its biological motivation, the main mathematical purpose of the present work is to rigorously confirm that such relaxing effects of
indirect chemoattractant production are actually not restricted to cases of chemotaxis-only systems, but rather
seem to form a much more general and robust feature of chemotactic interaction also in significantly more contexts
involving haptotactic migration mechanisms with all their potentially regularity-limiting properties due to
lack of haptoattractant diffusion.\abs
This will subsequently be examined in the framework of the problem
\be{0}
	\left\{ \begin{array}{ll}
	u_t=\Du \Delta u - \chi \nabla \cdot (u\nabla h) - \xi\nabla \cdot (u\nabla v) + f(u,v,w,h),
	\qquad & x\in\Omega, \ t>0, \\[1mm]
	h_t=\Dh \Delta h + g(u,v,w,h),
	\qquad & x\in\Omega, \ t>0, \\[1mm]
	v_t=-\alpha uv + v\phi(u,v,w,h) + \Phi(w),
	\qquad & x\in\Omega, \ t>0, \\[1mm]
	w_t=\beta u + w\psi(u,v,w,h),
	\qquad & x\in\Omega, \ t>0, \\[1mm]
	\Du \frac{\partial u}{\partial\nu} - \xi u\frac{\partial v}{\partial\nu}
	= \frac{\partial h}{\partial\nu}=0,
	\qquad & x\in\pO, \ t>0, \\[1mm]
	u(x,0)=u_0(x), \quad h(x,0)=h_0(x), \quad v(x,0)=v_0(x), \quad w(x,0)=w_0(x),
	\quad & x\in\Omega,
	\end{array} \right.
\ee
in a bounded domain $\Omega\subset \R^2$ with smooth boundary,
where the parameters $\chi,\xi, \Du$, $\Dh$, $\alpha$ and $\beta$ are assumed to be positive, 
and where for simplicity we suppose throughout this paper that with some $\vartheta\in (0,1)$,
\be{init}
	\left\{ \begin{array}{l}
	u_0 \in C^{2+\vartheta}(\bom)
	\quad \mbox{is nonnegative with } u_0\not\equiv 0, \\[1mm]
	h_0 \in C^{2+\vartheta}(\bom)
	\quad \mbox{is nonnegative}, \\[1mm]
	v_0 \in C^{2+\vartheta}(\bom)
	\quad \mbox{is positive in $\bom$ with $\frac{\partial v_0}{\partial\nu}=0$ on $\pO$, and that} \\[1mm]
	w_0 \in C^{2+\vartheta}(\bom)
	\quad \mbox{is positive in $\bom$.}
	\end{array} \right.
\ee
As for the parameter functions in (\ref{0}), in order to create a setup sufficiently general so as to include
both (\ref{01}) and (\ref{02}), we shall require that
\be{fgphipsi_reg}
	f,g,\phi \mbox{ and } \psi
	\quad \mbox{belong to} \quad 
	C^1([0,\infty)^4)
	\qquad \mbox{and} \qquad
	\Phi\in C^1([0,\infty)),
\ee
and are such that
\be{Hf}
	\left\{ \begin{array}{l}
	f_0(u) \le f(u,v,w,h) \le \Cf(v) \cdot (u+w+1)	
	\quad \mbox{for all } (u,v,w,h) \in [0,\infty)^4 \\[1mm]
	\mbox{with some $f_0\in C^1([0,\infty))$ such that $f_0(0)\ge 0$,
	and some nondecreasing $\Cf:[0,\infty)\to (0,\infty)$}
	\end{array} \right.
\ee
and
\be{Hg}
	\left\{ \begin{array}{l}
	|g(u,v,w,h)| \le \Cg(v) \cdot (w+h+1)
	\quad \mbox{for all } (u,v,w,h) \in [0,\infty)^4 \\[1mm]
	\mbox{with some nondecreasing $\Cg:[0,\infty)\to (0,\infty)$ \quad and} \\[1mm]
	g(0,v,w,0) \ge 0
	\quad \mbox{for all } (v,w) \in [0,\infty)^2, 
	\end{array} \right.
\ee
that
\be{Hphi}
	\left\{ \begin{array}{l}
	\phi(u,v,w,h) \le -\cphi \cdot w + \Cphi
	\quad \mbox{for all } (u,v,w,h) \in [0,\infty)^4, \\[1mm]
	|\phi_u(u,v,w,h)| \le \frac{\Cphi}{\sqrt{uv+1}}
	\quad \mbox{for all } (u,v,w,h) \in [0,\infty)^4, \\[1mm]
	|\phi_v(u,v,w,h)| \le \frac{\Cphi}{v+1} + \Cphi
	\quad \mbox{for all } (u,v,w,h) \in [0,\infty)^4, \\[1mm]
	|\phi_w(u,v,w,h)| \le \frac{\Cphi}{\sqrt{v+1}} + \Cphi
	\quad \mbox{for all } (u,v,w,h) \in [0,\infty)^4
	\quad \mbox{and} \quad \\[1mm]
	|\phi_h(u,v,w,h)| \le \frac{\Cphi}{\sqrt{v+1}} +\Cphi
	\quad \mbox{for all } (u,v,w,h) \in [0,\infty)^4 \\[1mm]
	\mbox{with some positive constants $\cphi$ and $\Cphi$,}
	\end{array} \right.
\ee
and that
\be{HPhi}
	\left\{ \begin{array}{l}
	0 \le \Phi(w) \le \CPhi
	\quad \mbox{for all } w\ge 0
	\quad \mbox{and} \\[1mm]
	w\Phi'^2(w) \le \CPhi \Phi(w)
	\quad \mbox{for all } w\ge 0,
	\end{array} \right.
\ee
and finally
\be{Hpsi}
	\left\{ \begin{array}{l}
	\psi(u,v,w,h) \le \Cpsi(v)
	\quad \mbox{for all } (u,v,w,h) \in [0,\infty)^4, \\[1mm]
	|\psi_u(u,v,w,h)| \le \frac{\Cpsi(v)}{\sqrt{uw+1}}
	\quad \mbox{for all } (u,v,w,h) \in [0,\infty)^4, \\[1mm]
	|\psi_v(u,v,w,h)| \le \Cpsi(v)
	\quad \mbox{for all } (u,v,w,h) \in [0,\infty)^4, \\[1mm]
	|\psi_w(u,v,w,h)| \le \frac{\Cpsi(v)}{w+1} 
	\quad \mbox{for all } (u,v,w,h) \in [0,\infty)^4
	\quad \mbox{and} \quad \\[1mm]
	|\psi_h(u,v,w,h)| \le \frac{\Cpsi(v)}{(w+1)^\gamma}
	\quad \mbox{for all } (u,v,w,h) \in [0,\infty)^4 \\[1mm]
	\mbox{with some nondecreasing $\Cpsi:[0,\infty)\to (0,\infty)$ and some $\gamma\in (0,\frac{1}{2})$.}
	\end{array} \right.
\ee
As can readily be verified, indeed both models (\ref{01}) and (\ref{02}) then become special cases of the PDE system
in (\ref{0}) whenever the parameters $\Du, \Dh, \chi, \xi, k_4$ and $k_6$ therein are positive, whereas 
$k_i$ for $i\in \{1,2,3,5,7,8\}$ and $\eta$ and, in particular, $\mu$ is merely required to be nonnegative.\abs
Our main results in this context then read as follows.
\begin{theo}\label{theo22}
  Let $\Omega\subset \R^2$ be a bounded domain with smooth boundary, assume that $\chi, \xi, \alpha, \beta, \Du$ and 
  $\Dh$ are positive, and suppose that $f$, $g$, $\phi$, $\Phi$ and $\psi$
  satisfy
  (\ref{fgphipsi_reg}), (\ref{Hf}), (\ref{Hg}), (\ref{Hphi}), (\ref{HPhi}) and (\ref{Hpsi}).
  Then for any choice of $(u_0,h_0,v_0,w_0)$ fulfilling (\ref{init}) with some $\vartheta\in (0,1)$, the problem 
  (\ref{0}) possesses a uniquely determined globally defined classical solution 
  $(u,h,v,w) \in (C^{2,1}(\bom\times [0,\infty))^4$ for which $u,h,v$ and $w$ are nonnegative.
\end{theo}
In particular, Theorem \ref{theo22} asserts that indeed no finite-time blow-up occurs in (\ref{0}) under
the above assumptions, and that in this regard the solution behavior in (\ref{0}) quite drastically differs from
that in the corresponding variants of (\ref{01}) and (\ref{02}) in which the second equation
is replaced with e.g.~$h_t=\Dh \Delta h + u$, and in which already in the semi-trivial case when $v\equiv 0$,
known results on the actually resulting two-component Keller-Segel system for $(u,h)$ assert finite-time blow-up
for some solutions (\cite{herrero_velazquez}).\abs
In line with this, all known results even on global existence, but also on qualitative properties,
in the original Chaplain-Lolas model (\cite{chaplain-lolas})
\be{03}
	\left\{ \begin{array}{l}
	u_t=\Delta u - \chi \nabla \cdot (u\nabla h) - \xi \nabla \cdot (u\nabla v) + \mu u(1-u-v), \\[1mm]
	h_t=\Delta h - h + u, \\[1mm]
	v_t=-hv + \eta v(1-u-v), 
	\end{array} \right.
\ee
seem to strongly rely on the assumption that $\mu$ be positive, thus guaranteeing the presence of a logistic-type quadratic
growth restriction on the cell density (\cite{cao_zamp2016}, \cite{pang_wang_m3as2018}, \cite{pang_wang_JDE},
\cite{tao_jmaa2009}, \cite{tao2014_arxiv}, \cite{taowin_SIMA2015}, \cite{wang_ke}).
Our results show that in the context of the variant (\ref{02}) of (\ref{03}) involving indirect chemoattractant
production, no such additional dampening
is necessary: Indeed, even when tissue remodeling is included by supposing that $\eta>0$ in (\ref{02}),
Theorem \ref{theo22} asserts global classical solvability in (\ref{02}) for actually any nonnegative value of $\mu$.\abs
The paper is structured as follows: 
After stating a result on local existence and extensibility as well as some preliminary estimates 
in Section \ref{sec:preliminaries}, in Section \ref{sec:ineq} we shall construct a quasi-energy functional for
(\ref{0}) and draw some immediate conclusions concerning regularity of solutions. 
The accordingly obtained estimates are used as a starting point for a Moser-type iterative argument yielding
$L^\infty$ bounds for the key solution component $u$ in Section \ref{sec:L-bounds}, and hence paving an essential
part of the way toward our proof of Theorem \ref{theo22} in Section \ref{sec:global-ex}. 
Finally, in Section \ref{sec:discussion} we illustrate the theoretical findings by numerical simulations 
of \eqref{02} and the corresponding model with direct signal production and provide some comments about the 
obtained results.\abs
\mysection{Local existence and basic estimates}\label{sec:preliminaries}
Following several precedents in the literature (\cite{fontelos_friedman_hu}, \cite{friedman_tello}, \cite{pang_wang_JDE}), 
in order to establish a preliminary result on local existence, but also to prepare our subsequent estimation
procedure, we note that on substituting
\be{z}
	z:=u e^{-\lambda v}
	\qquad \mbox{with} \qquad
	\lambda:=\frac{\xi}{\Du}>0,
\ee
the problem (\ref{0}) is equivalently transformed to 
\be{0z}
	\left\{ \begin{array}{ll}
	z_t =\Du e^{-\lambda v} \nabla \cdot (e^{\lambda v}\nabla z)
	- \chi e^{-\lambda v} \nabla \cdot (ze^{\lambda v} \nabla h) + e^{-\lambda v} f(z e^{\lambda v},v,w,h),
	\qquad & x\in \Omega, \ t>0, \\[1mm]
	h_t=\Dh \Delta h + g(ze^{\lambda v},v,w,h),
	\qquad & x\in\Omega, \ t>0, \\[1mm]
	v_t=-\alpha ve^{\lambda v}z + v\phi(ze^{\lambda v},v,w,h) + \Phi(w),
	\qquad & x\in\Omega, \ t>0, \\[1mm]
	w_t=\beta e^{\lambda v}z + w\psi(ze^{\lambda v},v,w,h),
	\qquad & x\in\Omega, \ t>0, \\[1mm]
	\frac{\partial z}{\partial\nu}
	= \frac{\partial h}{\partial\nu}=0,
	\qquad & x\in\pO, \ t>0, \\[1mm]
	z(x,0)=u_0(x)e^{\lambda v_0(x)}, \ h(x,0)=h_0(x), \ v(x,0)=v_0(x), \ w(x,0)=w_0(x),
	\ & x\in\Omega,
	\end{array} \right.
\ee
In this formulation, with respect to the construction of local-in-time solutions the problem (\ref{0}) indeed
becomes accessible to appropriate fixed point frameworks; 
by straightforward and minor adaptations of the corresponding arguments detailed e.g.~in \cite{pang_wang_JDE}, 
it is thereby possible to establish the following basic statement on unique solvability and extensibility.
\begin{lem}\label{lem_loc}
  Let $\chi, \xi, \alpha, \beta, \Du$ and 
  $\Dh$ be positive, let $f$, $g$, $\phi$, $\Phi$ and $\psi$
  comply with
  (\ref{fgphipsi_reg}), (\ref{Hf}), (\ref{Hg}), (\ref{Hphi}), (\ref{HPhi}) and (\ref{Hpsi}),
  and suppose that $(u_0,h_0,v_0,w_0)$ satisfies (\ref{init}) with some $\vartheta\in (0,1)$. 
  Then there exist $\tm\in (0,\infty]$ and a uniquely determined classical solution
  $(z,h,v,w) \in (C^{2,1}(\bom\times [0,\infty)))^4$ of (\ref{0z}) such that
  \be{ext}
	\mbox{either $\tm=\infty$, \quad or} \quad
	\limsup_{t\nearrow \tm} \Big\{ \|z(\cdot,t)\|_{L^\infty(\Omega)} + \|v(\cdot,t)\|_{W^{1,5}(\Omega)} \Big\}
	=\infty.
  \ee
  Moreover, we have $z>0,h\ge 0, v>0$ and $w>0$ in $\bom\times (0,\tm)$.
\end{lem}
Without any further explicit mentioning, throughout the sequel we shall suppose that the assumptions of Lemma \ref{lem_loc}
are satisfied, and that $(z,v,w,h)$ and $\tm\in (0,\infty]$ are as provided by the latter. 
Moreover, we shall tactitly switch between these variables and the quadruple $(u,v,w,h)$ solving
(\ref{0}) classically in $\Omega\times (0,\tm)$, as thereby defined through (\ref{z}).\abs
A first boundedness property of this solution is immediate.
\begin{lem}\label{lem1}
  The solution of (\ref{0}) satisfies
  \be{1.1}
	v(x,t) \le \Big\{ \|v_0\|_{L^\infty(\Omega)}  + \frac{\Cphi}{\CPhi}\Big\} \cdot e^{\Cphi t}
	\qquad \mbox{for all $x\in\Omega$ and } t\in (0,\tm).
  \ee
\end{lem}
\proof
  As from (\ref{0}), (\ref{Hphi}) and (\ref{HPhi}) we know that
  \bas
	v_t \le \Cphi v + \CPhi
	\qquad \mbox{in } \Omega\times (0,\tm),
  \eas
  by means of a simple comparison argument we conclude that
  \bas
	\|v(\cdot,t)\|_{L^\infty(\Omega)}
	&\le& \|v_0\|_{L^\infty(\Omega)} e^{\Cphi t}
	+ \CPhi \int_0^t e^{\Cphi (t-s)} ds
	\qquad \mbox{for all } t\in (0,\tm).
  \eas
  Herein estimating $\int_0^t e^{\Cphi (t-s)} ds \le \frac{1}{\Cphi} e^{\Cphi t}$ for $t\ge 0$, from this we readily obtain 
  (\ref{1.1}).
\qed
Next, thanks to (\ref{Hf}) and (\ref{Hg}) the first and fourth solution components in (\ref{0}) can at least controlled
with respect to their norm in $L^1(\Omega)$ in a fairly simple manner.
\begin{lem}\label{lem8}
  Let $T>0$. Then there exists $C(T)>0$ such that
  \be{8.1}
	\io u(\cdot,t) \le C(T)
	\quad \mbox{and} \quad
	\io w(\cdot,t) \le C(T)
	\qquad \mbox{for all }  t\in (0,\hatt)
  \ee
as well as
  \be{8.2}
	\io h(\cdot,t) \le C(T)
	\qquad \mbox{for all }  t\in (0,\hatt),
  \ee
where $\hatt:=\min\{T,\tm\}$.
\end{lem}
\proof
  Using that $v\le c_1(T):=\Big\{\|v_0\|_{L^\infty(\Omega)} +\frac{\CPhi}{\Cphi}\Big\} \cdot e^{\Cphi T}$ 
  in $\Omega\times (0,\hatt)$ due to Lemma \ref{lem1},
  after spatial integration in the first and the fourth equation in (\ref{0}) and adding the respective results 
  we may rely on (\ref{Hf}) and (\ref{Hpsi}) in estimating
  \bas
	\frac{d}{dt} \bigg\{ \io u + \io w \bigg\}
	&=& \io f(u,v,w,h) + \beta \io u + \io w\psi(u,v,w,h) \\
	&\le& \Cf(c_1(T)) \cdot \bigg\{ \io u + \io w + |\Omega| \bigg\}
	+ \beta \io u + \Cpsi(c_1(T)) \io w \\
	&\le& \Big\{ \Cf(c_1(T)) + \beta + \Cpsi(c_1(T)) \Big\} \cdot \bigg\{ \io u + \io w \bigg\}
	+ |\Omega| \Cf(c_1(T))
  \eas
  for $t\in (0,\hatt)$.
  A time integration of this linear ODI for $\io u + \io w$ directly yields (\ref{8.1}).
  Similarly, (\ref{Hg}) implies that
  \bas
	\frac{d}{dt} \io h
	= \io g(u,v,w,h)
	\le \Cg(c_1(T)) \cdot \bigg\{ \io w+ \io h + |\Omega| \bigg\}
  \eas
  for all $t\in (0,\hatt)$, so that (\ref{8.2}) becomes a consequence of (\ref{8.1}).
\qed
\mysection{A quasi-energy inequality}\label{sec:ineq}
The purpose of this section consists in the construction of an Lyapunov-like functional which through a corresponding
energy-dissipation inequality will provide some fundamental 
regularity information that will form the starting point of a series of {\em a priori} estimates which
in the presently considered spatially two-dimensional setting will finally allow for the conclusion that 
$(u,v,w,h)$ is actually global in time.\abs
As our first step in this direction, let us perform a standard testing procedure by which the crucial haptotactic
contribution to the first equation in (\ref{0}) is reduced to an $L^2$ inner product of gradients:
\begin{lem}\label{lem5}
  Let $\eta>0$. Then
  \bea{5.1}
	\frac{d}{dt} \io u\ln u
	+ \Du \io \frac{|\nabla u|^2}{u}
	&\le& \xi \io \nabla u\cdot\nabla v
	+ \frac{\chi^2}{2\eta} \io |\Delta h|^2 \nn\\
	& & + \eta \io u^2
	+ \Cf(\|v\|_{L^\infty(\Omega)}) \cdot \io u\ln u 
	+ 2\Cf(\|v\|_{L^\infty(\Omega)}) \cdot \io u \nn\\
	& & + \frac{\Cf^2(\|v\|_{L^\infty(\Omega)})}{2\eta} \cdot \io w^2
	+ \Cf(\|v\|_{L^\infty(\Omega)}) \cdot \io w \nn\\[2mm]
	& & + |\Omega| \Cf(\|v\|_{L^\infty(\Omega)})
	+ |\Omega| \cdot \|f_0\cdot\ln \|_{L^\infty((0,1))}
  \eea
  for all $t\in (0,\tm)$.
\end{lem}
\proof
  In the identity
  \bea{5.2}
	\frac{d}{dt} \io u\ln u 
	+ \Du \io \frac{|\nabla u|^2}{u}
	= \chi \io \nabla u\cdot\nabla h
	+ \xi \io \nabla u\cdot\nabla v
	+ \io f(u,v,w,h)\ln u
	+ \io f(u,v,w,h),
  \eea
  valid for all $t\in (0,\tm)$ due to (\ref{0}), we use Young's inequality to estimate
  \be{5.3}
	\chi \io \nabla u\cdot \nabla h
	= - \chi \io u\Delta h
	\le \frac{\eta}{2} \io u^2
	+ \frac{\chi^2}{2\eta} \io |\Delta h|^2
	\qquad \mbox{for all } t\in (0,\tm).
  \ee
  Moreover, by means of (\ref{Hf}) we see that for all $t\in (0,\tm)$,
  \be{5.5}
	\io f(u,v,w,h) \le \Cf(\|v\|_{L^\infty(\Omega)}) \cdot \bigg\{ \io u + \io w + |\Omega| \bigg\}
  \ee
  and
  \bas
	\io f(u,v,w,h) \ln u
	&\le& \Cf(\|v\|_{L^\infty(\Omega)}) \cdot \int_{\{u\ge 1\}} (u+w+1) \cdot \ln u
	+ \int_{\{u<1\}} f_0(u) \cdot \ln u \nn\\
	&\le& \Cf(\|v\|_{L^\infty(\Omega)}) \cdot \bigg\{ \io u\ln u + \io uw + \io u \bigg\}
	+ |\Omega| \cdot \|f_0\cdot\ln\|_{L^\infty((0,1))} \nn\\
	&\le& \Cf(\|v\|_{L^\infty(\Omega)}) \cdot \bigg\{ \io u\ln u + \io u\bigg\}
	+ \frac{\eta}{2} \io u^2
	+ \frac{\Cf^2(\|v\|_{L^\infty(\Omega)})}{2\eta} \io w^2 \nn\\[2mm]
	& & + |\Omega| \cdot \|f_0\cdot\ln\|_{L^\infty((0,1))},
  \eas
  because $\ln y \le y$ for all $y\ge 1$.
  In conjunction with (\ref{5.2})-(\ref{5.5}), this establishes (\ref{5.1}).
\qed
Now in order to achieve a precise cancelation of the integral in (\ref{5.1}) stemming from haptotactic 
cross-diffusion, inspired by several precedent works on similar types of interaction 
(see e.g.~\cite{corrias_perthame_zaag}, \cite{taowin_JDE2014} and \cite{SSW}), we track the evolution of the Dirichlet
integral associated with $\sqrt{v}$.
\begin{lem}\label{lem4}
  The solution of (\ref{0}) has the property that
  \bea{4.1}
	& & \hspace*{-30mm}
	\frac{d}{dt} \io \frac{|\nabla v|^2}{v}
	+ \frac{\cphi}{2} \io \frac{w}{v}|\nabla v|^2 \nn\\
	&\le& -2\alpha \io \nabla u\cdot \nabla v
	+ \frac{\alpha \Du}{2\xi} \io \frac{|\nabla u|^2}{u} \nn\\
	& & + \bigg\{ 2\Cphi \cdot (2\|v\|_{L^\infty(\Omega)}+1) + \frac{2\xi \Cphi^2}{\alpha \Du} \bigg\} 
		\cdot \io \frac{|\nabla v|^2}{v} \nn\\
	& & + \Big\{ \frac{4\Cphi^2 \cdot (\|v\|_{L^\infty(\Omega)}+1)}{\cphi} + \CPhi \Big\}
		\cdot \io \frac{|\nabla w|^2}{w} \nn\\
	& & + 2\Cphi \cdot (\|v\|_{L^\infty(\Omega)}+1) \cdot \io |\nabla h|^2
	\qquad \mbox{for all } t\in (0,\tm).
  \eea
\end{lem}
\proof
  According to the third equation in (\ref{0}), we have
  \bea{4.2}
	& & \hspace*{-10mm}
	\frac{d}{dt} \io \frac{|\nabla v|^2}{v} \nn\\
	&=& 2 \io \frac{\nabla v}{v} \cdot \nabla \Big\{ -\alpha uv + v\phi(u,v,w,h) + \Phi(w)\Big\}
	- \io \frac{|\nabla v|^2}{v^2} \cdot \Big\{-\alpha uv + v\phi(u,v,w,h)+\Phi(w)\Big\} \nn\\
	&=& -2\alpha \io \nabla u\cdot\nabla v
	- \alpha \io \frac{u}{v} |\nabla v|^2
	+ \io \frac{\phi(u,v,w,h)}{v} |\nabla v|^2 \nn\\
	& & + 2\io \phi_u(u,v,w,h)\nabla u\cdot\nabla v
	+ 2\io \phi_v(u,v,w,h) |\nabla v|^2 \\
	& & + 2 \io \phi_w(u,v,w,h) \nabla v\cdot \nabla w 
	+ 2\io \phi_h(u,v,w,h) \nabla v\cdot\nabla h \nn\\
	& & - \io \frac{\Phi(w)}{v^2} |\nabla v|^2
	+ 2\io \frac{\Phi'(w)}{v} \nabla v\cdot \nabla w
	\quad \mbox{for all } t\in (0,\tm).
  \eea
  Here from (\ref{Hphi}) we know that
  \be{4.3}
	\io \frac{\phi(u,v,w,h)}{v} |\nabla v|^2
	\le -\cphi \cdot \io \frac{w}{v} |\nabla v|^2
	+ \Cphi \cdot \io \frac{|\nabla v|^2}{v}
	\qquad \mbox{for all } t\in (0,\tm),
  \ee
  and that due to Young's inequality, for all $t\in (0,\tm)$ we have
  \bea{4.4}
	2 \io \phi_u(u,v,w,h) \nabla u \cdot\nabla v
	&\le& \frac{\alpha \Du}{2\xi} \io \frac{|\nabla u|^2}{u}
	+ \frac{2\xi}{\alpha\Du} \io u\phi_u^2(u,v,w,h) |\nabla v|^2 \nn\\
	&\le& \frac{\alpha \Du}{2\xi} \io \frac{|\nabla u|^2}{u}
	+ \frac{2\xi \Cphi^2}{\alpha\Du} \io \frac{|\nabla v|^2}{v}
  \eea
  and
  \bea{4.5}
	2 \io \phi_v(u,v,w,h)  |\nabla v|^2
	&\le& 2\Cphi \io \frac{|\nabla v|^2}{v}
	+ 2\Cphi \io |\nabla v|^2 \nn\\
	&\le& 2\Cphi \io \frac{|\nabla v|^2}{v}
	+ 2\Cphi \cdot \|v\|_{L^\infty(\Omega)} \io \frac{|\nabla v|^2}{v}
  \eea
  as well as
  \bea{4.6}
	2\io \phi_w(u,v,w,h) \nabla v\cdot \nabla w
	&\le& \frac{\cphi}{2} \io \frac{w}{v} |\nabla v|^2
	+ \frac{2}{\cphi} \io \frac{v\phi_w^2(u,v,w,h)}{w} |\nabla w|^2 \nn\\
	&\le& \frac{\cphi}{2} \io \frac{w}{v} |\nabla v|^2
	+ \frac{4\Cphi^2}{\cphi} \io \frac{|\nabla w|^2}{w}
	+ \frac{4\Cphi^2}{\cphi} \io \frac{v}{w} |\nabla w|^2 \nn\\
	&\le& \frac{\cphi}{2} \io \frac{w}{v} |\nabla v|^2
	+ \frac{4\Cphi^2 \cdot (\|v\|_{L^\infty(\Omega)}+1)}{\cphi} \io \frac{|\nabla w|^2}{w}
  \eea
  and
  \bea{4.7}
	2\io \phi_h(u,v,w,h) \nabla v\cdot\nabla h 
	&\le& \Cphi \io \frac{|\nabla v|^2}{v}
	+ \frac{1}{\Cphi} \io v\phi_h^2(u,v,w,h)|\nabla h|^2 \nn\\
	&\le& \Cphi \io \frac{|\nabla v|^2}{v}
	+ 2\Cphi \io |\nabla h|^2 
	+ 2\Cphi \io v|\nabla h|^2 \nn\\
	&\le& \Cphi \io \frac{|\nabla v|^2}{v}
	+ 2\Cphi\cdot (\|v\|_{L^\infty(\Omega)}+1) \cdot \io |\nabla h|^2.
  \eea
  Moreover, combining Young's inequality with (\ref{HPhi}) we can estimate
  \bea{4.8}
	- \io \frac{\Phi(w)}{v^2} |\nabla v|^2
	+ 2\io \frac{\Phi'(w)}{v} \nabla v\cdot \nabla w
	&\le& \io \frac{\Phi'^2(w)}{\Phi(w)} |\nabla w|^2 \nn\\
	&\le& \CPhi \io \frac{|\nabla w|^2}{w}
	\qquad \mbox{for all } \ t\in (0,\tm).
  \eea
  As the second summand on the right of (\ref{4.2}) is nonpositive, collecting (\ref{4.3})-(\ref{4.8}) we thus infer
  (\ref{4.1}) from (\ref{4.2}).
\qed
Next, several expressions on the right-hand sides of (\ref{5.1}) and (\ref{4.1}) need to be controlled in modulus.
Here the second integral on the right of (\ref{5.1}) can in fact be absorbed by the dissipation rate appearing in the
following inequality gained by means of a standard procedure.
\begin{lem}\label{lem6}
  For all $t\in (0,\tm)$,
  \bea{6.1}
	\frac{d}{dt} \io |\nabla h|^2 + \Dh \io |\Delta h|^2
	\le \frac{3\Cg^2(\|v\|_{L^\infty(\Omega)})}{\Dh} \cdot \bigg\{ \io w^2 + \io h^2 + |\Omega| \bigg\}.
  \eea
\end{lem}
\proof
  On testing the second equation in (\ref{0}) by $\Delta h$  and using Young's inequality in a standard manner, we obtain
  \bas
	\frac{1}{2} \frac{d}{dt} \io |\nabla h|^2 + \Dh \io |\Delta h|^2
	&=& -\io g(u,v,w,h) \Delta h \\
	&\le& \frac{\Dh}{2} \io |\Delta h|^2
	+ \frac{1}{2\Dh} \io g^2(u,v,w,h)
	\qquad \mbox{for all } t\in (0,\tm),
  \eas
  which implies (\ref{6.1}) due to the fact that by (\ref{Hg}) and again Young's inequality,
  \bas
	g^2(u,v,w,h)
	\le \Cg^2(\|v\|_{L^\infty(\Omega)}) \cdot (w+h+1)^2
	\le 3 \Cg^2(\|v\|_{L^\infty(\Omega)}) \cdot (w^2+h^2+1)
  \eas
  in $\Omega\times (0,\tm)$.
\qed
The second last summand in (\ref{4.1}), referring to the component $w$ with evolution governed by an ODE only, 
apparently cannot be expected to be absorbed by some suitable dissipation rate.
The following lemma indicates that at least some exponential control thereof will eventually
be possible.		
\begin{lem}\label{lem3}
  We have
  \bea{3.1}
	\frac{d}{dt} \io \frac{|\nabla w|^2}{w}
	&\le& (\beta+1) \io \frac{|\nabla u|^2}{u}
	+ \io \frac{w}{v} |\nabla v|^2 \nn\\
	& & + \Big\{ 3\Cpsi(\|v\|_{L^\infty(\Omega)}) + \Cpsi^2(\|v\|_{L^\infty(\Omega)})
	+ \|v\|_{L^\infty(\Omega)} \Cpsi^2(\|v\|_{L^\infty(\Omega)}) +1 \Big\} \cdot \io \frac{|\nabla w|^2}{w} \nn\\
	& & + \io w^2 
	+ \Cpsi^\frac{4}{1+2\gamma} (\|v\|_{L^\infty(\Omega)}) \cdot \io |\nabla h|^\frac{4}{1+2\gamma}
	\qquad \mbox{for all } t\in (0,\tm).
  \eea
\end{lem}
\proof
  Using the fourth equation in (\ref{0}), we compute
  \bea{3.2}
	\hspace*{-10mm}
	\frac{d}{dt} \io \frac{|\nabla w|^2}{w}
	&=& 2\io \frac{\nabla w}{w} \cdot \nabla \Big\{ \beta u + w\psi(u,v,w,h)\Big\}
	- \io \frac{|\nabla w|^2}{w^2} \cdot \Big\{\beta u + w\psi(u,v,w,h)\Big\} \nn\\
	&=& 2\beta \io \frac{1}{w} \nabla u\cdot\nabla w
	+ \io \frac{\psi(u,v,w,h)}{w} |\nabla w|^2
	- \beta \io \frac{u}{w^2} |\nabla w|^2 \nn\\
	& & +2\io \psi_u(u,v,w,h) \nabla u \cdot \nabla w
	+ 2 \io \psi_v(u,v,w,h) \nabla v\cdot\nabla w \nn\\
	& & + 2\io \psi_w(u,v,w,h) |\nabla w|^2
	+ 2\io \psi_h(u,v,w,h) \nabla w\cdot\nabla h
	\qquad \mbox{for all } t\in (0,\tm),
  \eea
  where by Young's inequality,
  \be{3.3}
	2\beta \io \frac{1}{w} \nabla u\cdot\nabla w
	- \beta \io \frac{u}{w^2}|\nabla w|^2
	\le \beta \io \frac{|\nabla u|^2}{u}
	\qquad \mbox{for all } t\in (0,\tm),
  \ee
  and where by (\ref{Hpsi}),
  \be{3.4}
	\io \frac{\psi(u,v,w,h)}{w} |\nabla w|^2
	\le \Cpsi(\|v\|_{L^\infty(\Omega)}) \cdot \io \frac{|\nabla w|^2}{w}
	\qquad \mbox{for all } t\in (0,\tm).
  \ee
  Apart from that, combining (\ref{Hpsi}) with Young's inequality we see that
  \bea{3.5}
	2 \io \psi_u(u,v,w,h) \nabla u \cdot\nabla w
	&\le& \io \frac{|\nabla u|^2}{u}
	+ \io u\psi_u^2(u,v,w,h) |\nabla w|^2 \nn\\
	&\le& \io \frac{|\nabla u|^2}{u}
	+ \Cpsi^2(\|v\|_{L^\infty(\Omega)}) \cdot \io \frac{|\nabla w|^2}{w}
  \eea
  and
  \bea{3.6}
	2\io \psi_v(u,v,w,h) \nabla v\cdot \nabla w
	&\le& \io \frac{w}{v} |\nabla v|^2
	+ \io \frac{v}{w} \psi_v^2(u,v,w,h) |\nabla w|^2 \nn\\
	&\le& \io \frac{w}{v}|\nabla v|^2
	+ \|v\|_{L^\infty(\Omega)} \Cpsi^2(\|v\|_{L^\infty(\Omega)}) \cdot \io \frac{|\nabla w|^2}{w}
  \eea
  as well as
  \bea{3.7}
	2 \io \psi_w(u,v,w,h) |\nabla w|^2
	\le 2\Cpsi(\|v\|_{L^\infty(\Omega)}) \cdot \io \frac{|\nabla w|^2}{w} 
  \eea
  and
  \bea{3.8}
	2\io \psi_h(u,v,w,h)\nabla w\cdot\nabla h
	&\le& \io \frac{|\nabla w|^2}{w}
	+ \io w\psi_h^2(u,v,w,h) |\nabla h|^2 \nn\\
	&\le& \io \frac{|\nabla w|^2}{w}
	+ \Cpsi^2(\|v\|_{L^\infty(\Omega)}) \cdot \io w^{1-2\gamma} |\nabla h|^2 \nn\\
	&\le& \io \frac{|\nabla w|^2}{w}
	+ \io w^2 + \Cpsi^\frac{4}{1+2\gamma}(\|v\|_{L^\infty(\Omega)}) \cdot \io |\nabla h|^\frac{4}{1+2\gamma}
  \eea
  for all $t\in (0,\tm)$.
  Inserting (\ref{3.3})-(\ref{3.8}) into (\ref{3.2}) directly yields (\ref{3.1}).
\qed
A final minor ingredient to our quasi-energy inequality is addressed in the following.
\begin{lem}\label{lem7}
  If $\eta>0$, then
  \be{7.1}
	\frac{d}{dt} \io w^2
	\le \eta \io u^2 + \Big\{ \frac{\beta^2}{\eta} + 2\Cpsi(\|v\|_{L^\infty(\Omega)})\Big\} \cdot \io w^2
	\qquad \mbox{for all } t\in (0,\tm).
  \ee
\end{lem}
\proof
  As
  \bas
	\frac{1}{2} \frac{d}{dt} \io w^2 = \beta \io uw + \io w^2 \psi(u,v,w,h)
	\qquad \mbox{for all } t\in (0,\tm)
  \eas
  by (\ref{0}), this follows by observing that
  \bas
	\beta \io uw \le \frac{\eta}{2} \io u^2 + \frac{\beta^2}{2\eta} \io w^2
	\qquad \mbox{for all } t\in (0,\tm)
  \eas
  due to Young's inequality, and that
  \bas
	\io w^2 \psi(u,v,w,h) \le \Cpsi(\|v\|_{L^\infty(\Omega)}) \cdot \io w^2
	\qquad \mbox{for all  }t\in (0,\tm)
  \eas
  according to (\ref{Hpsi}).
\qed
As a last preparation, let us make use of appropriate parabolic regularization features to estimate
terms of the form appearing in the last integral from (\ref{3.1}).
\begin{lem}\label{lem9}
  Let $q\in [1,4)$. Then there exists $\theta=\theta(q)\in (0,1)$ with the property
  that for all $T>0$ one can find $C(q,T)>0$ such that
  \be{9.1}
	\io |\nabla h|^q \le C(q,t) \cdot \bigg\{ \io |\Delta h|^2 \bigg\}^\theta
	\qquad \mbox{for all } t\in (0,\hatt),
  \ee
  where again $\hatt:=\min\{T,\tm\}$.
\end{lem}
\proof
We first note that according to (\ref{Hg}), Lemma \ref{lem1} and Lemma \ref{lem8} 
  we can fix positive constants $c_1(T), c_2(T)$ and $c_3(T)$
  such that
  \be{9.2}
	|g(u,v,w,h)| \le c_1(T) \cdot (w+h+1)
	\qquad \mbox{in } \Omega\times (0,\hatt),
  \ee
  and that
  \be{9.3}
	\io w \le c_2(T)
	\quad \mbox{and} \quad
	\io h \le c_3(T)
	\qquad \mbox{for all } t\in (0,\tm),
  \ee
  and observe that since $q<4$ it is possible to choose $r\equiv r(q)\in (1,2)$ suitably close to $2$ satisfying $r>q-2$.
  In the Duhamel representation
  \bas
	\nabla h(\cdot,t)
	= \nabla e^{t \Dh\Delta} h_0
	+ \int_0^t \nabla e^{(t-s)\Dh \Delta} g(u(\cdot,s),h(\cdot,s),v(\cdot,s),w(\cdot,s))ds,
	\qquad t\in (0,\tm),
  \eas
  we may then use standard $L^p$-$L^q$ estimates for the Neumann heat semigroup $(e^{\sigma\Delta})_{\sigma\ge 0}$
  (\cite[Lemma 1.3]{win_JDE2010}) to find $c_4(T)>0$ such that for all $t\in (0,\hatt)$,
  \bas
	\|\nabla h(\cdot,t)\|_{L^r(\Omega)} \le c_4(T) \|h_0\|_{W^{1,\infty}(\Omega)}
	+ c_4(T) \int_0^t (t-s)^{-\frac{3}{2}+\frac{1}{r}} 
		\Big\| g(u(\cdot,s),h(\cdot,s),v(\cdot,s),w(\cdot,s))\Big\|_{L^1(\Omega)} ds.
  \eas
  As (\ref{9.2}) and (\ref{9.3}) warrant that
  \bas
	\Big\| g(u(\cdot,s),h(\cdot,s),v(\cdot,s),w(\cdot,s))\Big\|_{L^1(\Omega)}
	\le c_5(T):=c_1(T) \cdot \Big(c_2(T)+c_3(T)+|\Omega|\Big)
	\qquad \mbox{for all } s\in (0,\hatt),
  \eas
  this readily entails that
  \be{9.4}
	\|\nabla h(\cdot,t)\|_{L^r(\Omega)}
	\le c_6(T):=c_4(T) \|h_0\|_{W^{1,\infty}(\Omega)}
	+ \frac{c_4(T) c_5(T) \cdot T^{\frac{1}{r}-\frac{1}{2}}}{\frac{1}{r}-\frac{1}{2}}
	\qquad \mbox{for all } t\in (0,\hatt)
  \ee
  with $c_6(T)$ being finite due to our restriction $r<2$.
  Now since a combination of the Gagliardo-Nirenberg inequality with elliptic regularity theory yields $c_7>0$ such that
  \bas
	\|\nabla h(\cdot,t)\|_{L^q(\Omega)}^q
	\le c_7\|\Delta h(\cdot,t)\|_{L^2(\Omega)}^{q-r}  \|\nabla h(\cdot,t)\|_{L^r(\Omega)}^r
	\qquad \mbox{for all } t\in (0,\tm),
  \eas
  from (\ref{9.4}) we immediately obtain (\ref{9.1}) with $\theta\equiv \theta(q):=\frac{q-r}{2}$ fulfilling
  $\theta\in (0,1)$ according to the inequality $r>q-2$.
\qed
We can now proceed to our detection of an energy-like structure in (\ref{0}), as expressed in the
following lemma.
\begin{lem}\label{lem77}
  Let $T>0$. Then there exist $a=a(T)>0, b>0$ and $C=C(T)>0$ such that for
  \be{77.1}
	\F(t):=\io u(\cdot,t)\ln u(\cdot,t) + a\io |\nabla h(\cdot,t)|^2 
	+\frac{\xi}{2\alpha} \io \frac{|\nabla v(\cdot,t)|^2}{v(\cdot,t)}
	+ b \io \frac{|\nabla w(\cdot,t)|^2}{w(\cdot,t)}
	+ \io w^2(\cdot,t),
	\qquad t\in [0,\hatt),
  \ee
  and
  \be{77.2}
	\D(t):=\io \frac{|\nabla u(\cdot,t)|^2}{u(\cdot,t)}
	+ \io |\Delta h(\cdot,t)|^2,
	\qquad t\in (0,\hatt),
  \ee
  we have
  \be{77.3}
	\frac{d}{dt} \F(t) + \frac{1}{C(T)} \cdot \D(t) \le C(T) \cdot \F(t) + C(T)
	\qquad \mbox{for all } t\in (0,\hatt),
  \ee
  where $\hatt:=\min\{T,\tm\}$.
\end{lem}
\proof
  Thanks to Lemma \ref{lem8} and Lemma \ref{lem1}, we can fix positive constants $c_1(T), c_2(T)$ and $c_3(T)$ such that
  \be{77.4}
	\io u \le c_1(T)
	\quad \mbox{and} \quad
	\io w  \le c_2(T)
	\qquad \mbox{for all } t\in (0,\hatt)
  \ee
  as well as
  \be{77.5}
	\|v(\cdot,t)\|_{L^\infty(\Omega)} \le c_3(T)
	\qquad \mbox{for all } t\in (0,\hatt).
  \ee
  Then in accordance with the Gagliardo-Nirenberg inequality, choosing $c_4>0$ such that
  \be{77.6}
	\io \varphi^4 \le c_4 \cdot \bigg\{ \io |\nabla \varphi|^2 \bigg\} \cdot \bigg\{ \io \varphi^2 \bigg\}
	+ c_4 \cdot \bigg\{ \io \varphi^2 \bigg\}^2
	\qquad \mbox{for all } \varphi\in W^{1,2}(\Omega),
  \ee
  we define 
  \be{77.7}
	\eta\equiv \eta(T):=\frac{\Du}{2c_1(T) c_4}
  \ee
  and take $a\equiv a(T)>0$ large enough fulfilling
  \be{77.8}
	\frac{\chi^2}{2\eta} \le \frac{a\Dh}{4}.
  \ee
  Finally picking $b>0$ small such that both
  \be{77.9}
	b(\beta+1) \le \frac{\Du}{4}
  \ee
  and
  \be{77.10}
b \le \frac{\xi \cphi}{4\alpha}
  \ee
  hold, we let $\F$ and $\D$ be as determined through (\ref{77.1}) and (\ref{77.2}) and claim that then (\ref{77.3}) is valid
  with some suitably large $C(T)>0$.\\
  To this end, we first take an appropriate linear combination of the inequalities provided 
  by Lemma \ref{lem5}, Lemma \ref{lem6}, Lemma \ref{lem4}, Lemma \ref{lem3} and Lemma 
  \ref{lem7}, which when applied to our particular value
  of $\eta$ namely show that
  \begin{align*} 
	\frac{d}{dt} \F(t)
	\le &\quad \bigg\{ -\Du \io \frac{|\nabla u|^2}{u}
	+ \xi \io \nabla u\cdot\nabla v
	+ \frac{\chi^2}{2\eta} \io |\Delta h|^2 \nn\\
	&  +\eta\io u^2 
	+ \Cf(c_3(T)) \cdot \io u\ln u
	+ 2\Cf(c_3(T)) \cdot c_1(T) \nn\\
	&  + \frac{\Cf^2(c_3(T))}{2\eta} \cdot \io w^2
	+ \Cf(c_3(T)) \cdot c_1(T) \nn\\
	&  + |\Omega| \Cf(c_3(T)) + |\Omega| \cdot \|f_0\cdot\ln\|_{L^\infty((0,1))} \bigg\} \nn
	\end{align*}
	 \bea{77.11}
	& & + \bigg\{ - a\Dh \io |\Delta h|^2	
	+ \frac{3a\Cg^2(c_3(T))}{\Dh} \cdot \io w^2
	+ \frac{3a\Cg^2(c_3(T))}{\Dh} \cdot \io h^2
	+ \frac{3a|\Omega| \Cg^2(c_3(T))}{\Dh} \bigg\} \nn\\
	& & + \bigg\{ - \frac{\xi \cphi}{4\alpha} \io \frac{w}{v}|\nabla v|^2 - \xi \io \nabla u\cdot\nabla v \nn\\
	& & + \frac{\Du}{4} \io \frac{|\nabla u|^2}{u}
	+ \frac{\xi}{2\alpha} \cdot \Big\{ 2\Cphi \cdot (2c_3(T)+1) + \frac{2\xi \Cphi^2}{\alpha\Du}\Big\} 
		\cdot \io \frac{|\nabla v|^2}{v} \nn\\
	& & + \frac{\xi}{2\alpha} \cdot \Big\{ \frac{4\Cphi^2 (c_3(T)+1)}{\cphi} + \CPhi \Big\} \io \frac{|\nabla w|^2}{w}
	+ \frac{\xi \Cphi \cdot (c_3(T)+1)}{\alpha} \io |\nabla h|^2 \bigg\} \nn\\
	& & + \bigg\{ b(\beta+1) \io \frac{|\nabla u|^2}{u}
	+ b\io \frac{w}{v} |\nabla v|^2 \nn\\
	& & + b\cdot \Big\{ 3\Cpsi(c_3(T)) + \Cpsi^2(c_3(T)) + c_3(T) \Cpsi^2(c_3(T)) +1\Big\} 
		\cdot \io \frac{|\nabla w|^2}{w} \nn\\
	& & + b\io w^2 + b\Cpsi^\frac{4}{1+2\gamma}(c_3(T)) \cdot \io |\nabla h|^\frac{4}{1+2\gamma} \bigg\} \nn\\
	& & + \bigg\{ \eta \io u^2 
	+ \Big\{ \frac{\beta^2}{\eta} + 2\Cpsi(c_3(T)) \Big\} \cdot \io w^2 \bigg\} \nn\\[2mm]
	&=& \Big\{ - \frac{3\Du}{4} + b(\beta+1)\Big\} \cdot \io \frac{|\nabla u|^2}{u}
	+ 2\eta \io u^2 \nn\\
	& & + \Big\{ -a\Dh + \frac{\chi^2}{2\eta}\Big\}\cdot \io |\Delta h|^2 
	+ c_5(T) \io |\nabla h|^\frac{4}{1+2\gamma} \nn\\
	& & + \Big\{ -\frac{\xi \cphi}{4\alpha} +b\Big\} \cdot \io \frac{w}{v} |\nabla v|^2 \nn\\
	& & + c_6(T) \io u\ln u 
	+ c_7(T) \io \frac{|\nabla v|^2}{v}
	+ c_8(T) \io \frac{|\nabla w|^2}{w}
	+ c_9(T) \io w^2 \nn\\
	& & + c_{10}(T) \io |\nabla h|^2 + c_{11}(T) \io h^2 
	+ c_{12}(T)
	\qquad \mbox{for all } t\in (0,\hatt)
  \eea
  with obvious definitions of $c_i(T)$ for $i\in\{5,...,12\}$, where we have made use of a favorable cancelation in some 
  contributions stemming from the haptotactic interaction.\\
  Now employing (\ref{77.6}) we see that due to (\ref{77.4}) and (\ref{77.7}) we have
  \bea{77.12}
	2\eta \io u^2
	&\le& \frac{c_4 \eta}{2} \cdot \bigg\{ \io \frac{|\nabla u|^2}{u} \bigg\} \cdot \bigg\{ \io u \bigg\}
	+ 2 c_4 \eta \cdot \bigg\{ \io u \bigg\}^2 \nn\\
	&\le& \frac{c_1(T) c_4 \eta}{2} \io \frac{|\nabla u|^2}{u} + 2c_1^2(T) c_4 \eta \nn\\
	&=& \frac{\Du}{4} \io \frac{|\nabla u|^2}{u} + 2c_1^2(T) c_4 \eta
	\qquad \mbox{for all } t\in (0,\hatt),
  \eea
  whereas Lemma \ref{lem9} says that as $\frac{4}{1+2\gamma}<4$, there exist $\theta\in (0,1)$ and $c_{13}(T)>0$ such that
  \bas
	c_5(T) \io |\nabla h|^\frac{4}{1+2\gamma} 
	\le c_{13}(T) \cdot \bigg\{ \io |\Delta h|^2 \bigg\}^\theta
	\qquad \mbox{for all } t\in (0,\hatt),
  \eas
  whence by Young's inequality,
  \be{77.13}
	c_5(T) \io |\nabla h|^\frac{4}{1+2\gamma}
	\le \frac{a\Dh}{4} \io |\Delta h|^2
	+ c_{14}(T)
	\qquad \mbox{for all } t\in (0,\hatt)
  \ee
  with some $c_{14}(T)>0$.
  Since moreover the Poincar\'e inequality provides $c_{15}(T)>0$ satisfying
  \bas
	c_{11}(T) \io h^2 \le c_{15}(T) \io |\nabla h|^2 + c_{15}(T)
	\qquad \mbox{for all } t\in (0,\hatt),
  \eas
  thanks to the restrictions in (\ref{77.8}), (\ref{77.9}) and (\ref{77.10}) we conclude from (\ref{77.11}), (\ref{77.12})
  and (\ref{77.13}) that there exists $c_{16}(T)>c_6(T)$ such that
  \bas
	\frac{d}{dt} \F(t)
	&\le& - \frac{\Du}{4} \io \frac{|\nabla u|^2}{u}
	- \frac{a\Dh}{2} \io |\Delta h|^2 \nn\\
	& & + c_6(T) \io u\ln u 
	+ c_{16}(T) \cdot \bigg\{ a\io |\nabla h|^2 + \frac{\xi}{2\alpha} \io \frac{|\nabla v|^2}{v}
	+ b\io \frac{|\nabla w|^2}{w} + \io w^2 \bigg\}
	+ c_{16}(T)
  \eas
  for all $t\in (0,\hatt)$. 
  Using that $y\ln y \ge -\frac{1}{e}$ for all $y>0$ and thus $\io u\ln u \ge -\frac{|\Omega|}{e}$
  for all $t\in (0,\tm)$, from this we finally obtain that for all $t\in (0,\hatt)$,
  \bas
	\frac{d}{dt} \F(t) + \min \Big\{ \frac{\Du}{4} \, , \, \frac{a\Dh}{2}\Big\} \cdot \D(t)
	&\le& c_{16}(T) \cdot \F(t)
	- \Big(c_{16}(T)-c_6(T)\Big) \cdot \io u\ln u
	+ c_{16}(T) \\
	&\le& c_{16}(T) \cdot \F(t)
	+ \frac{|\Omega| \cdot (c_{16}(T)-c_6(T))}{e} + c_{16}(T),
  \eas
  as intended.
\qed
Upon integration, the latter implies several {\em a priori} estimates, significantly going beyond those from
Lemma \ref{lem1} and Lemma \ref{lem8}, among which we explicitly state only those three inequalities that will be referred to 
later on.

\begin{lem}\label{lem10}
  For all $T>0$ there exists $C(T)>0$ such that again writing $\hatt:=\min\{T,\tm\}$ we have
  \be{10.1}
	\io u(\cdot,t) |\ln u(\cdot,t)| \le C(T)
	\qquad \mbox{for all } t\in (0,\hatt)
  \ee
  and 
  \be{10.2}
	\io |\nabla h(\cdot,t)|^2 \le C(T)
	\qquad \mbox{for all } t\in (0,\hatt)
  \ee
  as well as
  \be{10.3}
	\int_0^{\hatt} \io |\Delta h|^2 \le C(T).
  \ee
\end{lem}
\proof
  Upon integrating (\ref{77.3}) in time, we can find $c_1(T)>0$ such that with $\F$ and $\D$ taken from (\ref{77.1})
  and (\ref{77.2}) we have
  \bas
	\F(t) \le c_1(T)
	\quad \mbox{for all } t\in (0,\hatt)
	\qquad \mbox{and} \qquad
	\int_0^{\hatt} \D(t) dt \le c_1(T).
  \eas
  Once more using that $y\ln y\ge -\frac{1}{e}$ for $y>0$, from this we readily obtain the claimed inequalities
  as particular consequences.
\qed
As an immediate consequence, let us add the following observation about regularity of $h$.
\begin{cor}\label{cor100}
  Let $p\ge 1$ and $T>0$. Then there exists $C(p,T)>0$ such that with $\hatt:=\min\{T,\tm\}$,
  \be{100.1}
	\io h^p(\cdot,t) \le C(p,T)
	\qquad \mbox{for all } t\in (0,\hatt).
  \ee
\end{cor}
\proof
  As $W^{1,2}(\Omega)\hra L^p(\Omega)$, combining (\ref{8.2}) with (\ref{10.2}) immediately yields (\ref{100.1}).
\qed
\mysection{$L^\infty$ bounds for $u$}\label{sec:L-bounds}
In this section we intend to make use of the information from Lemma \ref{lem10} in order to finally achieve an {\em a priori}
bound for the quantity $z$ from (\ref{z}), and hence for $u$, with respect to the norm in $L^\infty(\Omega)$.
Here a key role will be played by the following implication of the estimate (\ref{10.1}) on a Gagliardo-Nirenberg-type
interpolation, as expressed in the following.
\begin{lem}\label{lem13}
  Let $p>1$ and $T>0$. Then for all $\eta>0$ one can find $C(\eta,p,T)>0$ such that
  \be{13.1}
	\io z^{p+1}
	\le \eta \io z^{p-2} |\nabla z|^2 + C(\eta,p,T)
	\qquad \mbox{for all } t\in (0,\hatt)
  \ee
  with $\hatt:=\min\{T,\tm\}$.
\end{lem}
\proof
  This follows in a standard manner from the boundedness property of $z$ in $L\log L(\Omega)$ 
  as implied by Lemma \ref{lem10} and Lemma \ref{lem1},
  by means of a refined interpolation inequality of Gagliardo-Nirenberg type, the latter going back to \cite{biler_HN}
  and extended to a version applicable to the present setting in \cite[Lemma A.5]{taowin_JDE2014} (for details, see 
  e.g.~\cite[p.800]{taowin_JDE2014}).
\qed
In order to make appropriate use of this, let us perform another well-established testing procedure to (\ref{0z}),
a basic outcome of which is the following.
\begin{lem}\label{lem11}
  Let $T>0$. Then there exists $C(T)>0$ such that for all $p\ge 2$ and all $t\in (0,\hatt)$ with $\hatt:=\min\{T,\tm\}$,
  \bea{11.1}
	\frac{d}{dt} \io e^{\lambda v} z^p
	+ \frac{p(p-1)\Du}{2} \io z^{p-2} |\nabla z|^2
	&\le& C(T) \cdot p \cdot \bigg\{ \io z^p +\io w^p +1 \bigg\} \nn\\
	& & + C(T) \cdot p^2 \cdot \io z^p |\nabla h|^2.
  \eea
\end{lem}
\proof
  By means of (\ref{0z}), we obtain
  \bea{11.2}
	\frac{d}{dt} \io e^{\lambda v} z^p
	&=& p \io e^{\lambda v} z^{p-1} \cdot \Big\{ \Du e^{-\lambda v} \nabla \cdot (e^{\lambda v}\nabla z)
	- \chi e^{-\lambda v} \nabla \cdot (z e^{\lambda v} \nabla h)
	+ e^{-\lambda v} f(z e^{\lambda v},h,v,w) \Big\} \nn\\
	& & + \lambda \io e^{\lambda v} z^p \cdot \Big\{ -\alpha v e^{\lambda v} z + v\phi(z e^{\lambda v},h,v,w)
		+ \Phi(w) \Big\} \nn\\
	&=& - p(p-1)\Du \io e^{\lambda v} z^{p-2} |\nabla z|^2
	+ p(p-1)\chi \io e^{\lambda v} z^{p-1} \nabla z\cdot\nabla h \nn\\
	& & + p \io z^{p-1} f(z e^{\lambda v},h,v,w)
	- \alpha \io v e^{2\lambda v} z^{p+1} \nn\\
	& & + \lambda \io v e^{\lambda v} z^p \phi(z e^{\lambda v},h,v,w)
	+ \lambda \io v e^{\lambda v} z^p \Phi(w)
  \eea
  for all $t\in (0,\tm)$, where by Young's inequality,  
  \bea{11.3}
	p(p-1)\chi \io e^{\lambda v} z^{p-1} \nabla z\cdot\nabla h
	&\le& \frac{p(p-1)\Du}{2} \io e^{\lambda v} z^{p-2} |\nabla z|^2 
	+ \frac{p(p-1)\chi^2}{2\Du} \io e^{\lambda v} z^p |\nabla h|^2 \nn\\
	&\le& \frac{p(p-1)\Du}{2} \io e^{\lambda v} z^{p-2} |\nabla z|^2 
	+ \frac{p^2\chi^2}{2\Du} e^{\lambda c_1(T)} \io z^p |\nabla h|^2
  \eea
  for all $t\in (0,\hatt)$, 
  with $c_1(T):=\|v\|_{L^\infty(\Omega\times (0,\hatt))}$ being finite according to Lemma \ref{lem1}.
  Furthermore, (\ref{Hf}) and Young's inequality warrant that for all $t\in (0,\hatt)$,
  \bea{11.4}
	p \io z^{p-1} f(z e^{\lambda v},h,v,w)
	&\le& p\Cf(c_1(T)) \cdot \bigg\{ \io e^{\lambda v} z^p + \io z^{p-1} w + \io z^{p-1} \bigg\} \nn\\
	&\le& p\Cf(c_1(T)) \cdot \bigg\{ e^{\lambda c_1(T)} \io z^p + \io z^p + \io w^p + \io z^p + |\Omega| \bigg\}, 
  \eea
  while from (\ref{Hphi}) and (\ref{HPhi}) we know that
  \be{11.5}
	\lambda \io v e^{\lambda v} \phi(z e^{\lambda v},h,v,w)
	\le \lambda c_1(T) e^{\lambda c_1(T)} 		
	\qquad \mbox{for all } t\in (0,\hatt)
  \ee
and that
  \be{11.6}
	\lambda \io e^{\lambda v} z^p \Phi(w)
	\le \lambda e^{\lambda c_1(T)} \cdot \CPhi \io z^p
	\qquad \mbox{for all } t\in (0,\hatt).
  \ee
  As clearly $e^{\lambda v}\ge 1$, by nonnegativity of $\alpha$ we therefore infer (\ref{11.1}) from (\ref{11.2}) when
  combined with (\ref{11.3}), (\ref{11.4}), (\ref{11.5}) and (\ref{11.6}).
\qed
Here a suitable control of the crucial rightmost summand in (\ref{11.1}) will rely, besides on Lemma \ref{lem13}, also on the following statement which can be viewed as partially generalizing Lemma \ref{lem7}.
\begin{lem}\label{lem14}
  Let $p\ge 1$ and $T>0$. Then writing $\hatt:=\min\{T,\tm\}$ and $\tau:=\frac{1}{2}\hatt$, with some $C(p,T)>0$ we have
  \be{14.1}
	\io w^p(\cdot,t) \le C(p,T) \int_\tau^t \io z^p
	+ C(p,T)
	\qquad \mbox{for all } t\in (\tau,\hatt).
  \ee
\end{lem}
\proof
  Choosing $c_1(T)>0$ large enough fulfilling $v\le c_1(T)$ in $\Omega\times (0,\hat T)$ according to Lemma \ref{lem1}, equation \eqref{0}, and 
  due to Young's inequality and (\ref{Hpsi}) we can estimate
  \bas
	\frac{d}{dt} \io w^p
	&=& p\beta \io e^{\lambda v} w^{p-1} z + p \io w^p \psi( e^{\lambda v}z,h,v,w) \\
	&\le& p\beta e^{\lambda c_1(T)} \io w^{p-1} z
	+ p\Cpsi(c_1(T)) \cdot \io w^p \\
	&\le& p\beta e^{\lambda c_1(T)} \cdot \bigg\{ \io w^p + \io z^p \bigg\}
	+ p\Cpsi(c_1(T)) \cdot \io w^p 
	\qquad \mbox{for all } t\in (0,\hatt),
  \eas
  so that with some $c_2(p,T)>0$ we have
  \bas
	\frac{d}{dt} \io w^p
	\le c_2(p,T) \io w^p + c_2(T) \io z^p
	\qquad \mbox{for all } t\in (0,\hatt).
  \eas
  An integration thereof shows that
  \bas
	\io w^p(\cdot,t)
	&\le& e^{c_2(p,T) \cdot (t-\tau)} \cdot \io w^p(\cdot,\tau)
	+ c_2(T) \int_\tau^t e^{c_2(p,T)\cdot (t-s)} \cdot \io z^p(\cdot,s) ds \\
	&\le& e^{c_2(p,T)\cdot T} \cdot \io w^p(\cdot,\tau)
	+ c_2(T) e^{c_2(p,T)\cdot T} \cdot \int_\tau^t \io z^p
	\qquad \mbox{for all } t\in (\tau,\hatt)
  \eas
  and hence implies (\ref{14.1}). 
\qed
We can thereby use Lemma \ref{lem11} along with Lemma \ref{lem13} to derive the following $L^p$ estimate for $z$,
at this stage yet involving bounds possibly depending on the finite number $p\ge 2$.
\begin{lem}\label{lem12}
  Let $p\ge 2$ and $T>0$. Then there exists $C(p,T)>0$ such that
  \be{12.1}
	\io z^p(\cdot,t) \le C(p,T)
	\qquad \mbox{for all } t\in (0,\hatt),
  \ee
  where again $\hatt:=\min\{T,\tm\}$.
\end{lem}
\proof
  From Lemma \ref{lem11} we obtain $c_1(p)>0$ and $c_2(p,T)>0$ such that
  \bea{12.2}
	\frac{d}{dt} \io e^{\lambda v} z^p
	+ c_1(p) \io z^{p-2} |\nabla z|^2
	&\le& c_2(p,T) \io e^{\lambda v} z^p
	+ c_2(p,T) \cdot \bigg\{ \io w^p + 1 \bigg\} \nn\\
	& & + c_2(p,T) \io z^p |\nabla h|^2
	\qquad \mbox{for all } t\in (0,\hatt),
  \eea
  again because $e^{\lambda v} \ge 1$.
  Here since Lemma \ref{lem10} provides $c_3(T)>0$ such that
  \bas
	\io |\nabla h|^2 \le c_3(T)
	\qquad \mbox{for all } t\in (0,\hatt),
  \eas
  by means of Young's inequality and the Gagliardo-Nirenberg inequality we infer that with some $c_4(p,T)>0$ we have
  \bas
	c_2(p,T) \io z^p |\nabla h|^2
	&\le& c_2(p,T) \io z^{p+1}
	+ c_2(p,T) \io |\nabla h|^{2p+2} \\
	&\le& c_2(p,T) \io z^{p+1}
	+ c_4(p,T) \cdot \bigg\{ \io |\Delta h|^{p+1} \bigg\} \cdot \bigg\{ \io |\nabla h|^2 \bigg\}^\frac{p+1}{2} \\
	&\le& c_2(p,T) \io z^{p+1}
	+ (c_3(T))^\frac{p+1}{2} c_4(p,T) \io |\Delta h|^{p+1}
	\qquad \mbox{for all } t\in (0,\hatt).
  \eas
  Therefore, (\ref{12.2}) shows that there exists $c_5(p,T)>0$ such that $y(t):=\io e^{\lambda v(\cdot,t)} z^p(\cdot,t)$,
  $t\in [0,\tm)$, satisfies
  \be{12.3}
	y'(t) + c_1(p) \io z^{p-2} |\nabla z|^2
	\le c_2(p,T) y(t) + c_5(p,T) \cdot \bigg\{ \io z^{p+1} + \io w^{p+1} + \io |\Delta h|^{p+1} +1 \bigg\}
  \ee
  for all $t\in (0,\hatt)$ and thus, upon integration,
  \bea{12.4}
	& & \hspace*{-20mm}
	y(t) + c_1(p) \int_{\tau}^t \io z^{p-2} |\nabla z|^2 \nn\\
	&\le& y(t) + c_1(p) \int_{\tau}^t e^{c_2(p,T)(t-s)} \cdot \io z^{p-2}(\cdot,s) |\nabla z(\cdot,s)|^2 ds \nn\\
	&\le& y(\tau) e^{c_2(p,T)\cdot (t-\tau)} \nn\\
	& & + c_5(p,T) \int_{\tau}^t e^{c_2(p,T)\cdot (t-s)} \cdot 
	\bigg\{ \io z^{p+1}(\cdot,s) + \io w^{p+1}(\cdot,s) + \io |\Delta h(\cdot,s)|^{p+1} +1 \bigg\} ds \nn\\
	&\le& y(\tau) e^{c_2(p,T)\cdot T} \nn\\
	& & + c_5(p,T) e^{c_2(p,T)\cdot T} \cdot 
	\bigg\{ \int_{\tau}^t \io z^{p+1} + \int_{\tau}^t \io w^{p+1} + \int_\tau^t \io |\Delta h|^{p+1} + T \bigg\}
  \eea
  for all $t\in (\tau,\hatt)$, where again we have set $\tau:=\frac{1}{2}\hatt$.
  As thus $\tau$ is positive, a well-known result on maximal Sobolev regularity in the parabolic subproblem of (\ref{0})
  satisfied by $h$ (\cite{giga_sohr}) becomes applicable so as to yield $c_6(p,T)>0$ satisfying
  \bas
	\int_\tau^t \io |\Delta h|^{p+1} 
	&\le& c_6(p,T) \int_\tau^t \io |g(e^{\lambda v}z,h,v,w)|^{p+1} + c_6(p,T) \\
	&\le& c_6(p,T) \cdot (p+1) \Cg^{p+1}(\|v\|_{L^\infty(\Omega\times (0,\hatt))}) \cdot
	\int_\tau^t \io (w^{p+1} + h^{p+1} +1 )
	+ c_6(p,T)
  \eas
  for all $t\in (\tau,\hatt)$ because of (\ref{Hg}). 
  Since Corollary \ref{cor100} and Lemma \ref{lem14} provide $c_7(p,T)>0$ and $c_8(p,T)>0$ such that
  \bas
	\int_\tau^t \io h^{p+1} \le c_7(p,T)
	\qquad \mbox{for all } t\in (\tau,\hatt)
  \eas
  and that
  \bas
	\int_\tau^t \io w^{p+1} \le c_8(p,T) \cdot \bigg\{ \int_\tau^t \io z^{p+1} +1 \bigg\}
	\qquad \mbox{for all } t\in (\tau,\hatt),
  \eas
  this means that with some $c_9(p,T)>0$ we have
  \bas
	\int_\tau^T \io |\Delta h|^{p+1} \le c_9(p,T) \cdot \bigg\{ \int_\tau^t \io z^{p+1} + 1 \bigg\}
	\qquad \mbox{for all } t\in (\tau,\hatt),
  \eas
  so that from (\ref{12.4}) we infer the existence of $c_{10}(p,T)>0$ satisfying
  \be{12.7}
	y(t) + c_1(p) \int_\tau^t \io z^{p-2} |\nabla z|^2
	\le c_{10}(p,T) \int_\tau^t \io z^{p+1}
	+ c_{10}(p,T) 
	\qquad \mbox{for all } t\in (\tau,\hatt).
  \ee
  We now employ Lemma \ref{lem13} to see that with some $c_{11}(p,T)>0$,
  \bas
	c_{10}(p,T) \int_\tau^t \io z^{p+1} \le c_1(p) \int_\tau^t \io z^{p-2} |\nabla z|^2
	+ c_{11}(p,T)
	\qquad \mbox{for all } t\in (\tau,\hatt),
  \eas
  whence (\ref{12.7}) ensures that
  \bas
	y(t) \le c_{10}(p,T) + c_{11}(p,T)
	\qquad \mbox{for all } t\in (\tau,\hatt)
  \eas
  and that thus, clearly, (\ref{12.1}) holds.
\qed
By adapting a well-established Moser-type iteration (\cite{alikakos}, \cite{taowin_subcrit}) to the present context, however,
one can readily turn the latter into estimates in $L^\infty$.
\begin{lem}\label{lem15}
  Given $T>0$, one can find $C(T)>0$ such that with $\hatt:=\min\{T,\tm\}$ we have
  \be{15.1}
	\|z(\cdot,t)\|_{L^\infty(\Omega)} \le C(T)
	\qquad \mbox{for all } t\in (0,\hatt).
  \ee
\end{lem}
\proof
  We first fix any $p_\star>2$ and then obtain upon combining Lemma \ref{lem12} with 
  Lemma \ref{lem14}, Lemma \ref{lem1}, (\ref{Hg}) and Corollary \ref{cor100} that
  \bas
	\int_\tau^t \io |g|^{p_\star} \le c_1(T)
	\qquad \mbox{for all } t\in (\tau,\hatt)
  \eas
  with $\tau:=\frac{1}{2}\hatt$ and some $c_1(T)>0$.
  As a consequence thereof, standard regularization features of the Neumann heat semigroup (\cite[Lemma 1.3]{win_JDE2010},
  \cite[Lemma 4.1]{horstmann_win}) entail boundedness of $\nabla h$ in $\Omega\times (\tau,\hatt)$,
  so that from Lemma \ref{lem11} we infer the existence of $c_2(T)>0$ and $c_3(T)>0$ such that for $p:=2^{k+1}$
  and any nonnegative integer $k$,
  \bas
	\frac{d}{dt} \io e^{\lambda v} z^p
	+ c_2(T) \io |\nabla z^\frac{p}{2}|^2
	\le c_3(T) \cdot p^2 \cdot \bigg\{ \io z^p + \io w^p +1 \bigg\}
	\qquad \mbox{for all } t\in (\tau,\hatt).
  \eas
  On integrating and recalling Lemma \ref{lem1} and Lemma \ref{lem14}, we see that with some $c_4(T)>0$, for any such $p$
  and arbitrary $t\in (\tau,\hatt)$ we have
  \bea{15.2}	
	\io z^p(\cdot,t) + c_2(T) \int_\tau^t \io |\nabla z^\frac{p}{2}|^2
	&\le& \io e^{\lambda v(\cdot,t)} z^p(\cdot,t)
	+ c_2(T) \int_\tau^t \io |\nabla z^\frac{p}{2}|^2 \nn\\
	&\le& \io e^{\lambda v(\cdot,\tau)} z^p(\cdot,\tau) \nn\\
	& & + c_3(T) \cdot p^2 \cdot \int_\tau^t \io z^p
	+ c_3(T) \cdot p^2 \cdot \int_\tau^t \io w^p + c_3(T) \cdot p^2 \nn\\
	&\le& c_4(T) \io z^p(\cdot,\tau)
	+ c_4(T)\cdot p^2 \cdot \int_\tau^t \io z^p	
	+ c_3(T) \cdot p^2.
  \eea
  The remaining part now follows a well-established reasoning: 
  According to the Gagliardo-Nirenberg inequality and Young's inequality,
  we can find $c_5(T)>0$ such that introducing the numbers
  \bas
	M_k:=\max \bigg\{ 1 \, , \, \sup_{t\in (\tau,\hatt)} \io z^{p_k}(\cdot,t) \bigg\},
	\qquad k\in \{0,1,2,...\},
  \eas
  all finite due to Lemma \ref{lem12}, for $p=p_k$ and each $k\in \{1,2,3,...\}$ we have
  \bas
	c_4(T) \cdot p^2 \cdot \int_\tau^t \io z^p
	&\le& c_5(T) \cdot p^2 \cdot 
	\int_\tau^t \|\nabla z^\frac{p}{2}(\cdot,s)\|_{L^2(\Omega)} \|z^\frac{p}{2}(\cdot,s)\|_{L^1(\Omega)} ds \\
	&\le& c_5(T) \cdot p^2 \cdot M_{k-1} \cdot \int_\tau^t \|\nabla z^\frac{p}{2}(\cdot,s)\|_{L^2(\Omega)} ds \\
	&\le& c_2(T) \int_\tau^t \io |\nabla z^\frac{p}{2}|^2
	+ c_5^2(T) T \cdot p^4 \cdot M_{k-1}^2 
	\qquad \mbox{for all } t\in (\tau,\hatt),
  \eas
  whence (\ref{15.2}) entails that for some $c_6(T)>0$,
  \bas
	M_k 
	\le c_4(T) \io z^{p_k}(\cdot,\tau)
	+ c_6(T) \cdot p_k^4 \cdot M_{k-1}^2
	\qquad \mbox{for all } k\ge 1.
  \eas
  By means of a standard recursive argument, both when $p_k^4 M_{k-1}^2 \le \|z(\cdot,\tau)\|_{L^\infty(\Omega)}^{p_k}$
  for infinitely many $k\ge 1$, and as well in the opposite case this can readily be seen to imply the existence of
  $c_7(T)>0$ such that
  \bas
	\|z(\cdot,t)\|_{L^\infty(\Omega)} \le c_7(T)
	\qquad \mbox{for all } t\in (\tau,\hatt),
  \eas
  from which the claim immediately follows.
\qed
\mysection{Global extensibility. Proof of Theorem \ref{theo22}}\label{sec:global-ex}
Having thus ruled out blow-up of the first quantity appearing in the second alternative from (\ref{ext}),
it hence remains to derive appropriate bounds for the haptotactic gradient.
A first observation relates the latter to some spatio-temporal regularity properties of $z$ and $h$.
\begin{lem}\label{lem17}
  Let $T>0$ and $q\ge 1$. Then there exists $C(q,T)>0$ fulfilling
  \be{17.1}
	\|\nabla v(\cdot,t)\|_{L^q(\Omega)} 
	+ \|\nabla w(\cdot,t)\|_{L^q(\Omega)}
	\le C(q,T) \int_\tau^T \Big\{ \|\nabla z(\cdot,s)\|_{L^q(\Omega)}
	+ \|\nabla h(\cdot,s)\|_{L^q(\Omega)} \Big\} ds
	+ C(q,T)
  \ee
  for all $t\in (\tau,\hatt)$, with $\hatt:=\min\{T,\tm\}$ and $\tau:=\frac{1}{2}\hatt$.
\end{lem}
\proof
  Differentiating in (\ref{0}) and recalling (\ref{z}), we see that throughout $\Omega\times (0,\tm)$,
  \be{17.2}
	\partial_t \nabla v 
	= (-\alpha v + v\phi_u e^{\lambda v}) \nabla z
	+ (-\alpha \lambda uv - \alpha u + \lambda uv \phi_u + v\phi_v+\phi) \nabla v
	+ (v\phi_w +\Phi'(w)) \nabla w + v\phi_h \nabla h
  \ee
  and
  \be{17.3}
	\partial_t \nabla w
	= (\beta e^{\lambda v} + we^{\lambda v} \psi_u)\nabla z 
	+ (\beta\lambda u + \lambda uw + w\psi_v)\nabla v
	= (w\psi_w + \psi)\nabla w	
	+ w\phi_h \nabla h,
  \ee
  where we have suppressed the argument $(u,v,w,h)$ in $\phi$, $\psi$ and the derivatives thereof.
  Now as a consequence of Lemma \ref{lem1}, Lemma \ref{lem15} and (\ref{z}), our requirements (\ref{Hphi}), 
  (\ref{HPhi}) and (\ref{Hpsi})
  guarantee that herein all the functions 
  $-\alpha v + v\phi_u e^{\lambda v}$, $-\alpha \lambda uv - \alpha u + \lambda uv \phi_u + v\phi_v+\phi$,
  $v\phi_w+\Phi'(w)$, $v\phi_h$, $\beta e^{\lambda v} + we^{\lambda v} \psi_u$,
  $\beta\lambda u + \lambda uw + w\psi_v$, $w\psi_w + \psi$ and $w\psi_h$ are bounded in $\Omega\times (0,\hatt)$,
  so that (\ref{17.2}) and (\ref{17.3}) imply that with some $c_1(T)>0$ we have
  \bas
	\|\nabla v(\cdot,t)\|_{L^q(\Omega)}
	&=& \bigg\| \nabla v(\cdot,\tau)
	+ \int_\tau^t \partial_t \nabla v(\cdot,s) ds \bigg\|_{L^q(\Omega)} \\
	&\le& c_1(T) + c_1(T) \int_\tau^t \Big\{ \|\nabla z(\cdot,s)\|_{L^q(\Omega)}
	+ \|\nabla v(\cdot,s)\|_{L^q(\Omega)}\\
	&&\qquad \qquad \qquad + \|\nabla w(\cdot,s)\|_{L^q(\Omega)}+ \|\nabla h(\cdot,s)\|_{L^q(\Omega)} \Big\} ds
  \eas
  and, similarly, 
  \bas
	\|\nabla w(\cdot,t)\|_{L^q(\Omega)}
	\le c_1(T) + c_1(T) \int_\tau^t \Big\{ \|\nabla z(\cdot,s)\|_{L^q(\Omega)}
	+ \|\nabla v(\cdot,s)\|_{L^q(\Omega)}\\
		+ \|\nabla w(\cdot,s)\|_{L^q(\Omega)}
	+ \|\nabla h(\cdot,s)\|_{L^q(\Omega)} \Big\} ds
  \eas
  for all $t\in (\tau,\hatt)$. Adding these inequalities and invoking Gronwall's lemma readily leads to (\ref{17.1}).
\qed
Here the crucial ingredient containing $\nabla z$ will ultimately be controlled by using the following 
result which extends a corresponding statement from \cite[Lemma 3.14]{taowin_JDE2014} to the present more
complex system (\ref{0z}).
\begin{lem}\label{lem16}
  Let $T>0$. Then there exists $C(T)>0$ with the property that
  \be{16.1}
	\io |\nabla z(\cdot,t)|^2 \le C(T)
	\quad \mbox{for all } t\in (\tau,\hatt)
	\qquad \mbox{and} \qquad
	\int_\tau^T \io |\Delta z|^2 \le C(T),
  \ee
  where again $\hatt:=\min\{T,\tm\}$ and $\tau:=\frac{1}{2}\hatt$.
\end{lem}
\proof
  Defining
  \bas
	A(x,t):=\lambda \Du \nabla v - \chi\nabla h
	\qquad \mbox{and} \qquad
	B(x,t):=-\chi z\Delta h - \chi \lambda z \nabla v\cdot\nabla h
	+ e^{-\lambda v} f(e^{\lambda v}z,v,w,h)
  \eas
  for $(x,t)\in \Omega\times (0,\tm)$, we see that the first equation in (\ref{0z}) simplifies to the identity
  \bas
	z_t=\Du \Delta z + A(x,t)\cdot\nabla z + B(x,t),
	\qquad x\in\Omega, \ t\in (0,\tm),
  \eas
  which we multiply by $-\Delta z$ and integrate to find using Young's inequality that
  for all $t_0\in [0,\tm)$ and any $t\in (t_0,\tm)$,
  \bea{16.3}
	\io |\nabla z(\cdot,t)|^2 
	+ \Du \int_{t_0}^t \io |\Delta z|^2
	&=& \io |\nabla z(\cdot,t_0)|^2
	- \int_{t_0}^t \io (A\cdot\nabla z)\Delta z
	- \int_{t_0}^t \io B\Delta z \nn\\
	&\le& \io |\nabla z(\cdot,t_0)|^2
	+ \Du \int_{t_0}^t \io |\Delta z|^2
	+ \frac{1}{2\Du} \int_{t_0}^t \io |A\cdot\nabla z|^2 \nn\\
	& & + \frac{1}{2\Du} \int_{t_0}^t \io B^2.
  \eea
  Here we note that given $T>0$, in view of the boundedness property of $z$ in $\Omega\times (0,\hatt)$ asserted 
  by Lemma \ref{lem15} we can fix $c_1(T)>0$ such that due to Young's inequality,
  \bea{16.4}
	\frac{1}{\Du} \int_{t_0}^t \io B^2
	&\le& c_1(T) \int_{t_0}^{\hatt} \io |\Delta h|^2
	+ c_1(T) \int_{t_0}^t \io |\nabla v\cdot\nabla h|^2 
	+ c_1(T) \nn\\
	&\le& c_1(T) \int_{t_0}^{\hatt} \io |\Delta h|^2
	+ c_1(T) \int_{t_0}^t \io |\nabla v|^4
	+ \frac{c_1(T)}{4} \int_{t_0}^t \io |\nabla h|^4 \nn\\[2mm]
	& & + c_1(T)
	\qquad \qquad \mbox{for all $t_0\in [0,\hatt)$ and } t\in (t_0,\hatt).
  \eea
  Moreover, as the Gagliardo-Nirenberg inequality together with elliptic estimates and
  Lemma \ref{lem15} provides $c_2>0$ and $c_3(T)>0$ such that
  \bea{16.44}
	\io |\nabla z|^4
	&\le& c_2 \|\Delta z\|_{L^2(\Omega)}^2 \|z\|_{L^\infty(\Omega)}^2 
	+ c_2 \|z\|_{L^\infty(\Omega)}^4 \nn\\
	&\le& c_3(T) \|\Delta z\|_{L^2(\Omega)}^2
	+ c_3(T)
	\qquad \mbox{for all } t\in (0,\hatt),
  \eea
  by combining the Cauchy-Schwarz inequality with Young's inequality we obtain
  \bas
	\frac{1}{\Du} \int_{t_0}^t \io |A\cdot\nabla z|^2
	&\le& \frac{1}{\Du} \int_{t_0}^t \|A(\cdot,s)\|_{L^4(\Omega)}^2 \|\nabla z(\cdot,s)\|_{L^4(\Omega)}^2 ds \nn\\
	&\le& \frac{\sqrt{c_3(T)}}{\Du} \int_{t_0}^t \|A(\cdot,s)\|_{L^4(\Omega)}^2
	\sqrt{\|\Delta z(\cdot,s)\|_{L^2(\Omega)}^2 + 1 } ds \nn\\
	&\le& \frac{\Du}{4} \cdot \bigg\{ \int_{t_0}^t \io |\Delta z|^2 + (t-t_0)\bigg\}
	+ \frac{c_3(T)}{\Du^3} \int_{t_0}^t \io |A|^4 \nn\\
	&\le& \frac{\Du}{4} \int_{t_0}^t \io |\Delta z|^2
	+ \frac{\Du T}{4} \nn\\
	& & + 4c_3(T) \lambda^4 \Du \int_{t_0}^t \io |\nabla v|^4 \\
	& & + \frac{4c_3(T) \chi^4}{\Du^3} \int_{t_0}^t \io |\nabla h|^4 
	\qquad \mbox{for all $t_0\in [0,\hatt)$ and } t\in (t_0,\hatt),
  \eas
  which along with (\ref{16.4}) shows that (\ref{16.3}) implies the inequality
  \bea{16.5}
	& & \hspace*{-20mm}
	\io |\nabla z(\cdot,t)|^2	
	+ \frac{\Du}{4} \int_{t_0}^t \io |\Delta z|^2 \nn\\
	&\le& \io |\nabla z(\cdot,t_0)|^2 
	+ c_4(T) \int_\tau^{\hatt} \io |\Delta h|^2
	+ c_4(T) \int_\tau^{\hatt} \io |\nabla h|^4 \nn\\
	& & + c_4(T) \int_{t_0}^t \io |\nabla v|^4
	+ c_4(T)
	\qquad \mbox{for all $t_0\in [\tau,\hatt)$ and } t\in (t_0,\hatt)
  \eea
  with some appropriately large $c_4(T)>0$.\\
  To proceed from this, we observe that as a consequence of Lemma \ref{lem17} we can find $c_5(T)>0$ such that
  \bas
	\io |\nabla v|^4
	\le c_5(T) 
	+ c_5(T) \int_\tau^t \io |\nabla z|^4
	+ c_5(T) \int_\tau^{\hatt} \io |\nabla h|^4
	\qquad \mbox{for all } t\in (\tau,\hatt),
  \eas
  and that another application of the Gagliardo-Nirenberg inequality in conjunction with Lemma \ref{lem10} provides
  $c_6>0$ and $c_7(T)>0$ fulfilling
  \bas
	\io |\nabla h|^4
	&\le& c_6 \cdot \bigg\{ \io |\Delta h|^2 \bigg\} \cdot \bigg\{ \io |\nabla h|^2 \bigg\} \\
	&\le& c_7(T) \cdot \io |\Delta h|^2
	\qquad \mbox{for all } t\in (0,\hatt).
  \eas
  Therefore, again through Lemma \ref{lem10}, (\ref{16.5}) reduces to
  \bea{16.6}
	\io |\nabla z(\cdot,t)|^2
	+ \frac{\Du}{4} \int_{t_0}^t \io |\Delta z|^2
	&\le& \io |\nabla z(\cdot,t_0)|^2
	+ c_8(T) \int_{t_0}^t \int_\tau^s \io |\nabla z(x,\sigma)|^4 dxd\sigma ds \nn\\
	& & + c_8(T)
	\qquad \mbox{for all $t_0\in [\tau,\hatt)$ and } t\in (t_0,\hatt),
  \eea
  where by the Fubini theorem and, again, (\ref{16.44}),
  \bea{16.7}	
	& &  \hspace*{-18mm}
	c_8(T) \int_{t_0}^t \int_\tau^s \io |\nabla z(x,\sigma)|^4 dxd\sigma ds \nn\\
	&=& c_3(T) c_8(T) \int_{t_0}^t (t-\sigma)\cdot \io |\nabla z(x,\sigma)|^4 dxd\sigma  \nn\\
	& & + c_8(T) \cdot (t-t_0) \int_\tau^{t_0} \io |\nabla z|^4 \nn\\
	&\le& c_8(T) \cdot (t-t_0) \cdot \int_{t_0}^t \io |\nabla z|^4 \nn\\
	& & + c_8(T) \cdot (t-t_0) \cdot \int_\tau^{t_0} \io |\nabla z|^4 \nn\\
	&\le& c_3(T) c_8(T) \cdot (t-t_0) \cdot \int_{t_0}^t \io |\Delta z|^2 +c_8(T)C_3(T)T(t-t_0)\nn\\
	& & + c_3(T) c_8(T) T \int_\tau^{t_0} \io |\Delta z|^2 +c_8(T)C_3(T)T(t_0-\tau )
	\quad \mbox{for all $t_0\in [\tau,\hatt)$ and } t\in (t_0,\hatt).
  \eea
  We now let $\delta>0$ be small enough such that $c_3(T) c_8(T)\delta \le \frac{\Du}{8}$, and choose $N\in\N$
  and $(t_i)_{i\in\{1,...,N\}} \subset [\tau,\hatt]$ such that $t_1=\tau, t_N=\hatt$ and $0<t_{i+1}-t_i \le\delta$
  for all $i\in \{1,...,N-1\}$.
  Then inserting (\ref{16.7}) into (\ref{16.6}) shows that 
  \bas
	I_i:=\sup_{t\in (t_i,t_{i+1})} \io |\nabla z(\cdot,t)|^2
	\quad \mbox{and} \quad
	J_i:=\frac{\Du}{8} \int_{t_i}^{t_{i+1}} \io |\Delta z|^2,
	\qquad i\in \{1,...,N\},
  \eas
  along with $I_0:=\io |\nabla z(\cdot,\tau)|^2$ and $J_0:=0$ satisfy
  \bas
	\max \{I_i,J_i\}
	\le I_{i-1} + c_3(T) c_8(T) T \cdot \sum_{j=0}^{i-1} J_j + c_8(T)
	\qquad \mbox{for all } i\in \{1,...,N-1\},
  \eas
  whence if we let 
  \bas
	K_i:=\sum_{k=1}^i \max\{I_k,J_k\} +1
	\qquad \mbox{for } i\in\{1,...,N-1\},
  \eas
  then 
  \bas
	K_i
	&\le& \sum_{k=1}^i I_{k-1}
	+ c_3(T) c_8(T) T \cdot \sum_{k=1}^i \sum_{j=0}^{k-1} J_j + ic_8(T) +1 \\
	&\le& I_0 + \sum_{j=1}^{i-1} I_j
	+ c_3(T) c_8(T)T \cdot \sum_{j=1}^{i-1} \sum_{k=j+1}^i J_j
	+ (N-1) c_8(T) +1 \\
	&=& \sum_{j=1}^{i-1} \Big\{ I_j + c_3(T) c_8(T)T \cdot (i-j) J_j\Big\} + I_0 + (N-1)c_8(T) + 1 \\
	&\le& \Big\{ 1+ (N-2) c_3(T) c_8(T)T \Big\} \cdot K_{i-1}
	+ I_0 + (N-1) c_8(T) +1 \\[2mm]
	&\le& c_9(T) \cdot K_{i-1}
	\qquad \mbox{for all } i\in \{1,...,N-1\}
  \eas
  with $c_9(T):=\max\{ 1+(N-2)c_3(T) c_8(T)T \, , \, I_0 +(N-1) c_8(T) +1\}$.
  Therefore, $K_i \le K_0 \cdot c_9^{N-1}(T)$ for all $i\in \{1,...,N-1\}$, which 
  in view of the definitions of $(I_i)_{i\in\{1,...,N-1\}}$ and $(J_i)_{i\in\{1,...,N-1\}}$ yields (\ref{16.1}).
\qed
Combining this with Lemma \ref{lem17} readily implies the following.
\begin{lem}\label{lem18}
  For all $q\ge 1$ and $T>0$ one can fix $C(q,T)>0$ such that with $\hatt:=\min\{T,\tm\}$ and $\tau:=\frac{1}{2}\hatt$ we have
  \be{18.1}
	\|\nabla v(\cdot,t)\|_{L^q(\Omega)} 
	\le C(q,T)
	\qquad \mbox{for all } t\in (\tau,\hatt).
  \ee
\end{lem}
\proof
  From Lemma \ref{lem16} in conjunction with elliptic regularity theory we obtain $c_1(T)>0$ such that
  \bas
	\int_\tau^{\hatt} \|z(\cdot,t)\|_{W^{2,2}(\Omega)}^2 dt \le c_1(T),
  \eas
  which we combine with the continuity of the embedding $W^{2,2}(\Omega) \hra W^{1,q}(\Omega)$ and the Cauchy-Schwarz 
  inequality to conclude that with some $c_2(q)>0$ we have
  \bas
	\int_\tau^{\hatt} \|\nabla z(\cdot,t)\|_{L^q(\Omega)} dt
	&\le& c_2(q) \int_\tau^{\hatt} \|z(\cdot,t)\|_{W^{2,2}(\Omega)} dt \\
	&\le& c_2(q) \sqrt{T} \cdot 
	\bigg\{ \int_\tau^{\hatt} \|z(\cdot,t)\|_{W^{2,2}(\Omega)}^2 dt \bigg\}^\frac{1}{2} \\[2mm]
	&\le& \sqrt{c_1(T)} c_2(q) \cdot \sqrt{T}.
  \eas
  Since Lemma \ref{lem10} in quite a similar fashion yields the existence of $c_3(q,T)>0$ such that
  \bas
	\int_\tau^{\hatt} \|\nabla h(\cdot,t)\|_{L^q(\Omega)} \le c_3(q,T),
  \eas
  the claimed estimate is a consequence of Lemma \ref{lem17}.
\qed
Now our main result has actually already been established:\abs
\proofc of Theorem \ref{theo22}. \quad
  Thanks to Lemma \ref{lem_loc}, and in particular the extensibility criterion (\ref{ext}) therein,
  we only need to combine the outcome of Lemma \ref{lem15} with an application of Lemma \ref{lem18} to $q:=5$.
\qed
%
%
%


\section{Simulations and discussion}\label{sec:discussion}

\begin{figure}[h!]
		\includegraphics[width=0.23\textwidth]{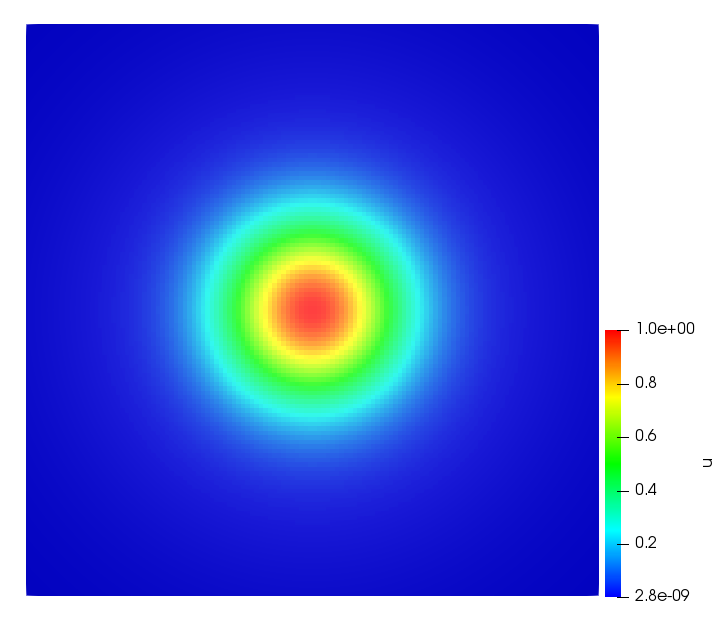}\quad 
		\includegraphics[width=0.23\textwidth]{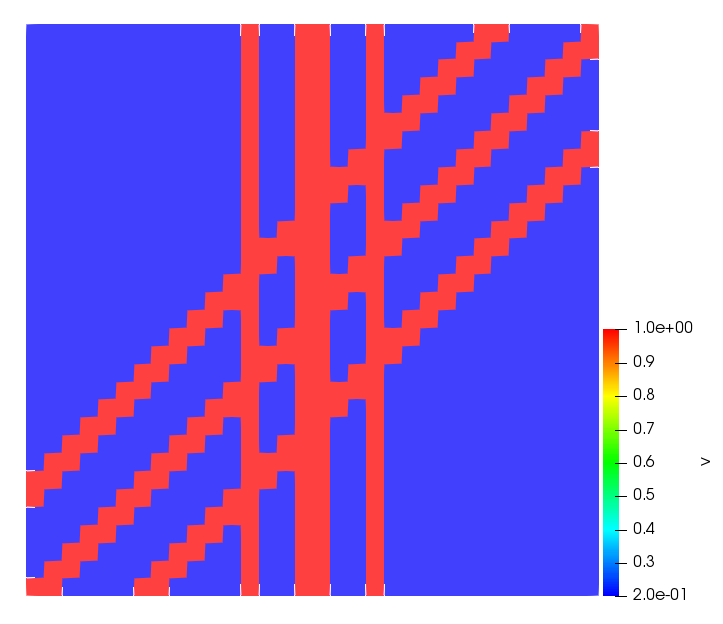}
		\quad \includegraphics[width=0.23\textwidth]{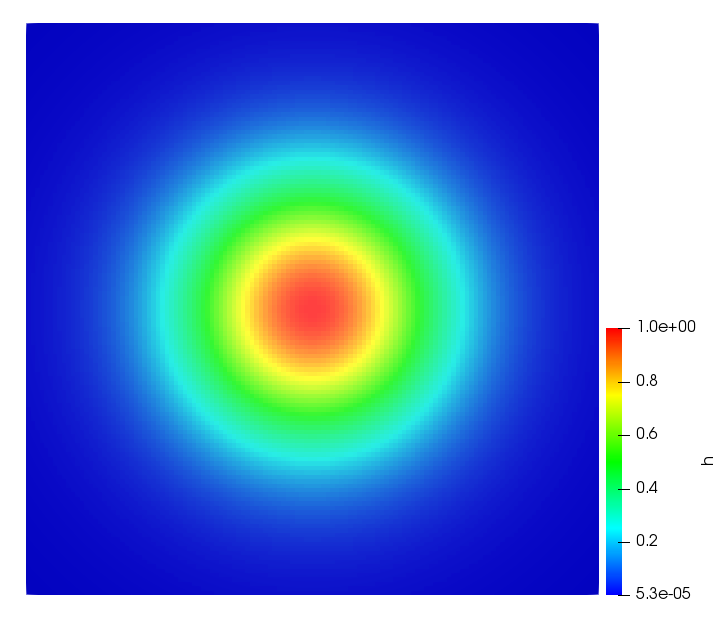}
		\quad \includegraphics[width=0.23\textwidth]{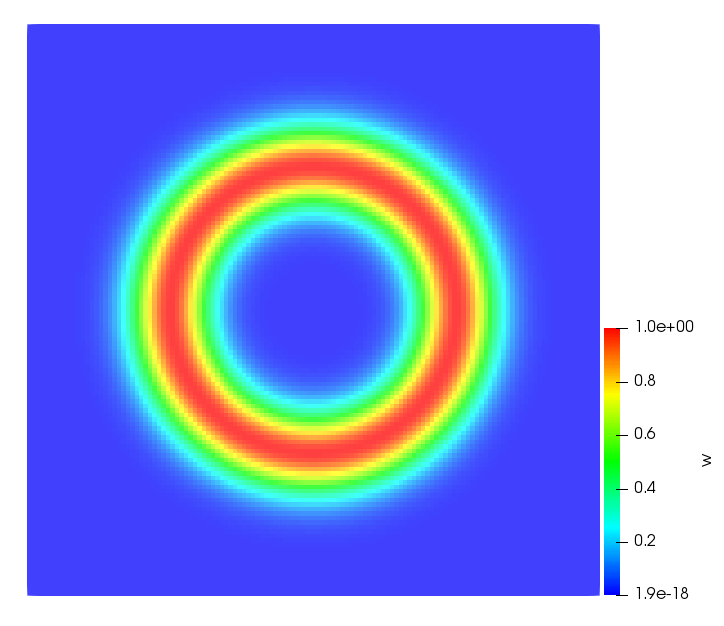}\\[1ex]
		\includegraphics[width=0.23\textwidth]{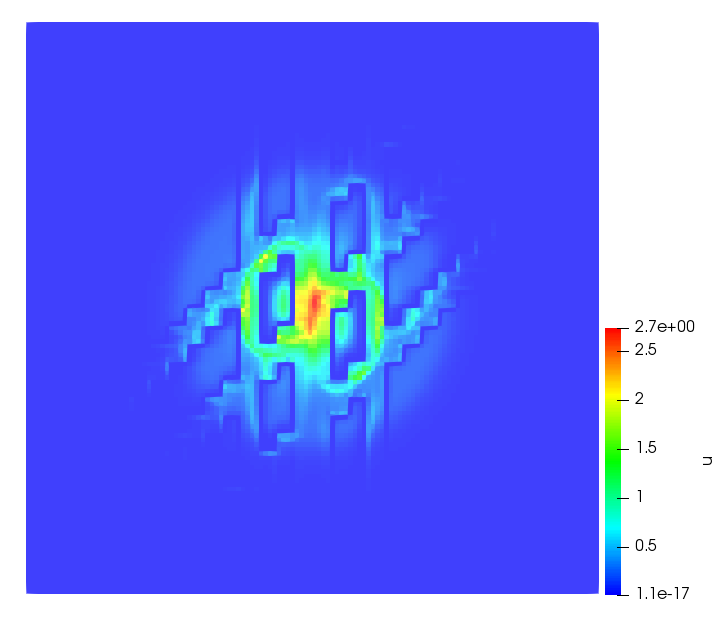}
		\quad \includegraphics[width=0.225\textwidth]{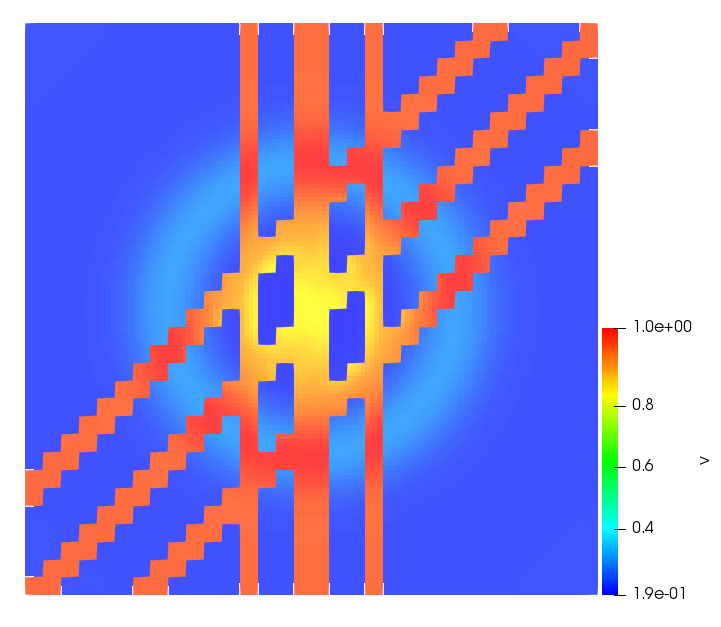}
		\quad \includegraphics[width=0.23\textwidth]{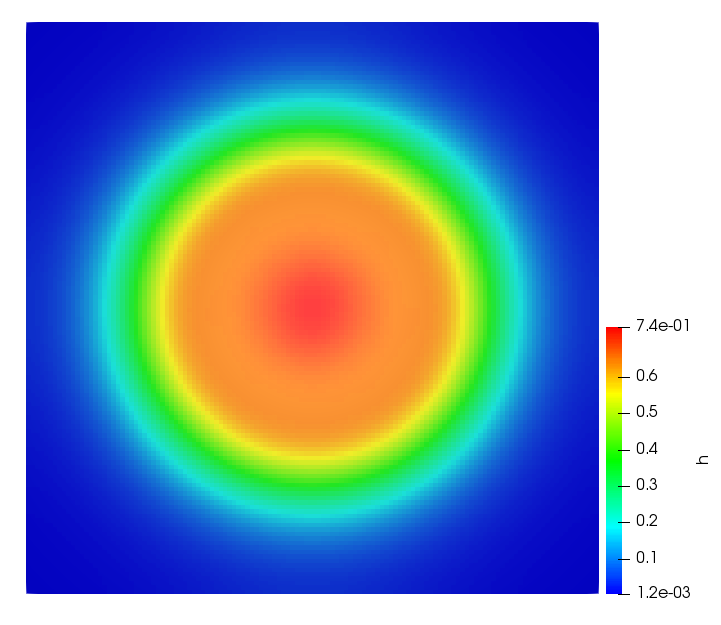}
		\quad \includegraphics[width=0.23\textwidth]{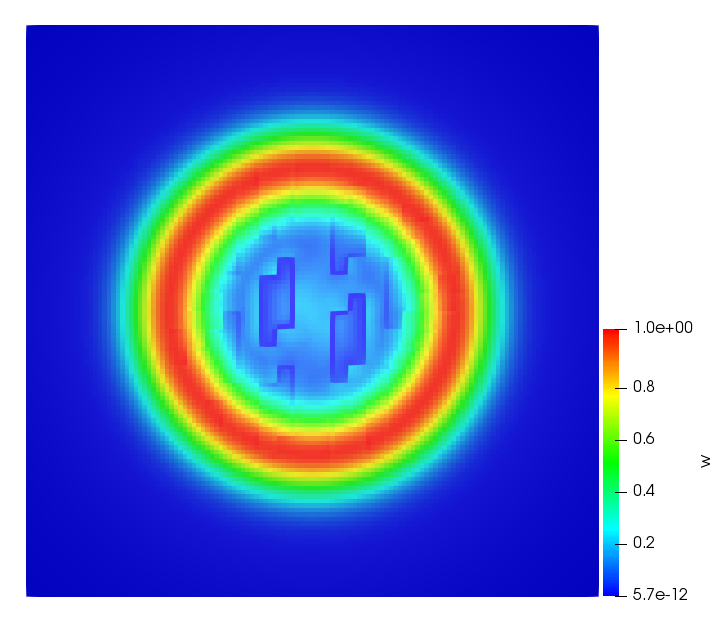}\\[1ex]
		\includegraphics[width=0.23\textwidth]{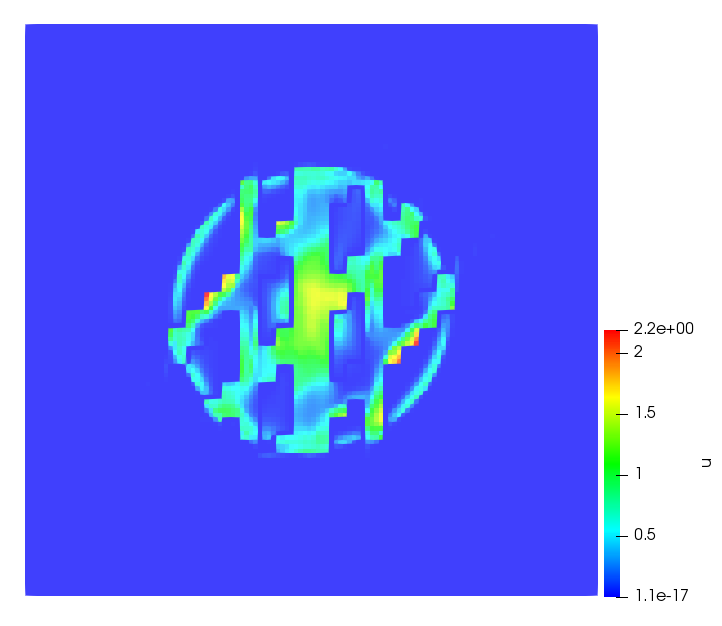}
		\quad \includegraphics[width=0.225\textwidth]{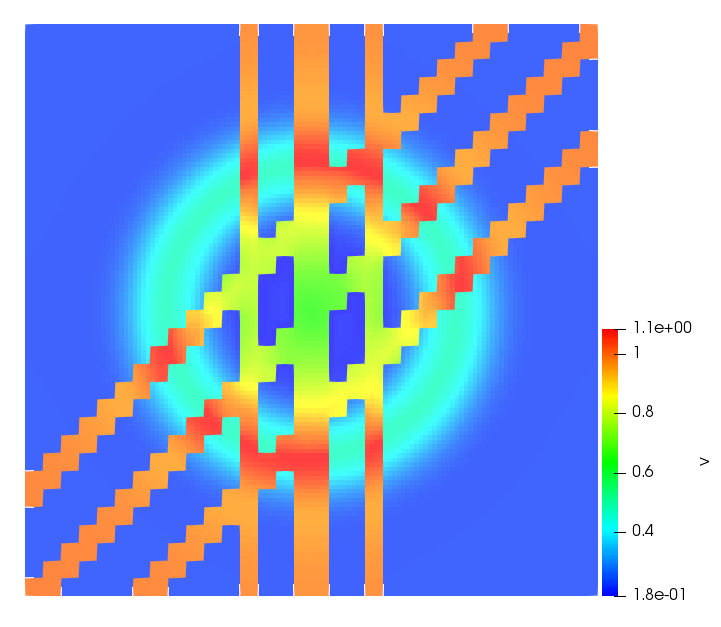}
		\quad \includegraphics[width=0.23\textwidth]{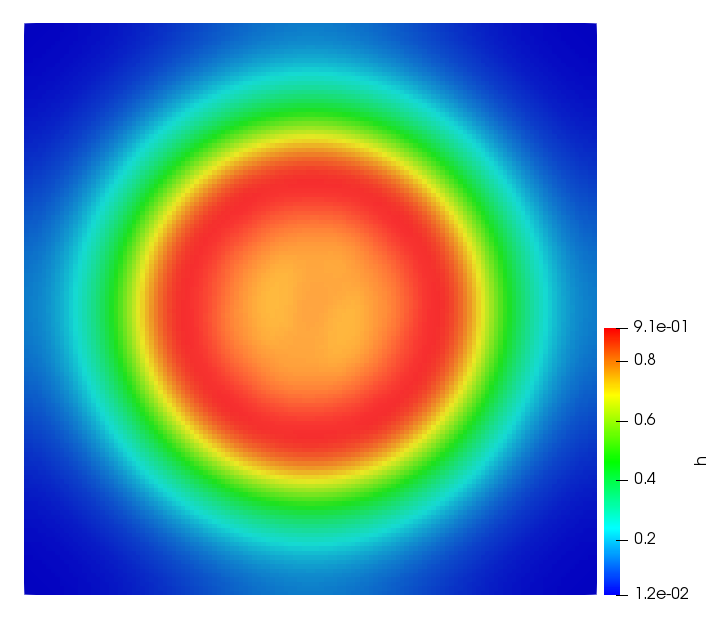}
		\quad \includegraphics[width=0.23\textwidth]{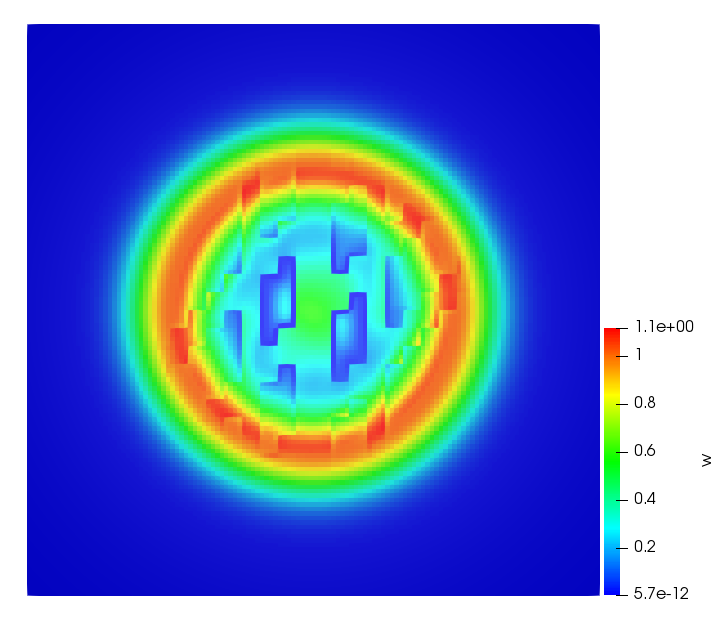}\\[1ex]
		\includegraphics[width=0.23\textwidth]{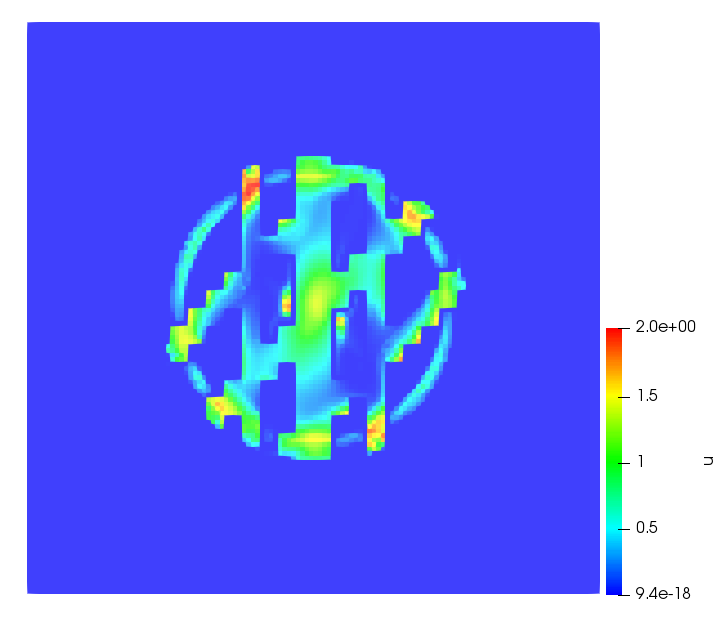}
		\quad \includegraphics[width=0.225\textwidth]{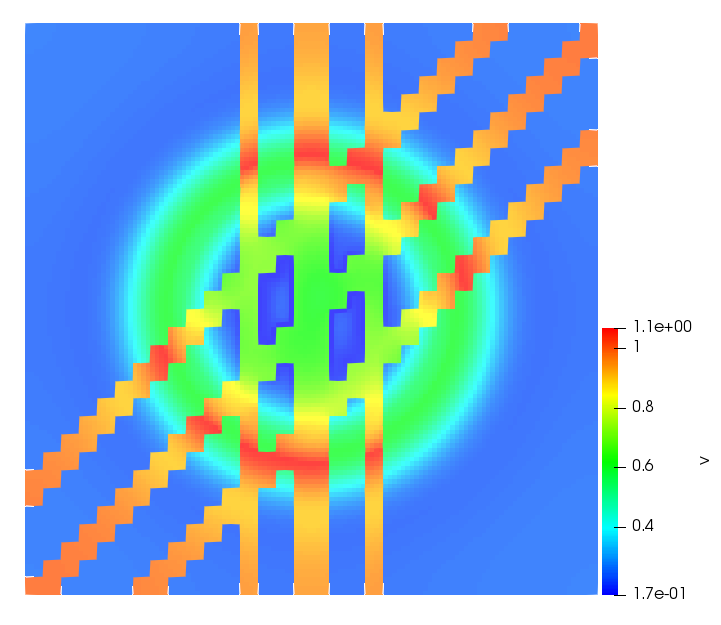}
		\quad \includegraphics[width=0.23\textwidth]{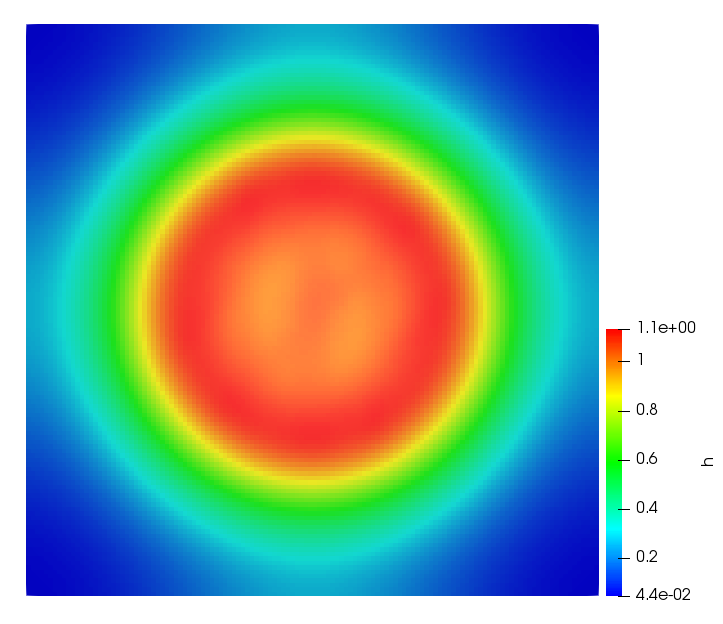}
		\quad \includegraphics[width=0.23\textwidth]{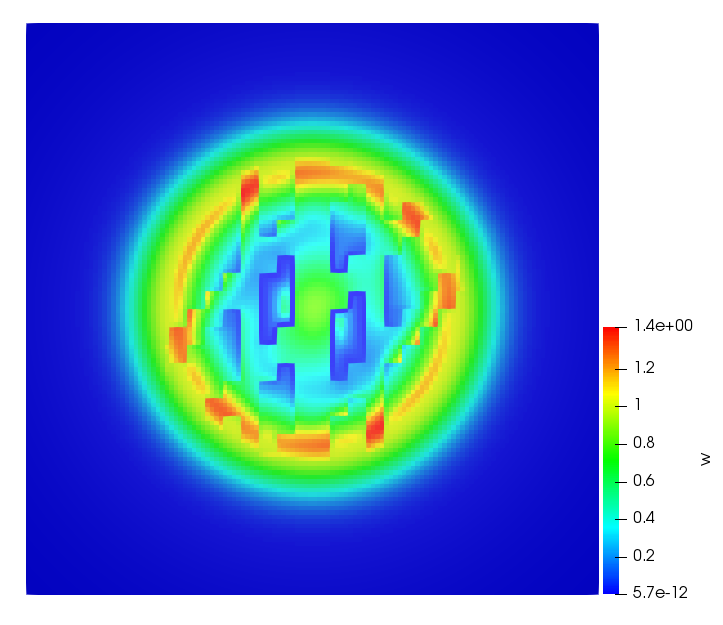}\\[1ex]
		\includegraphics[width=0.23\textwidth]{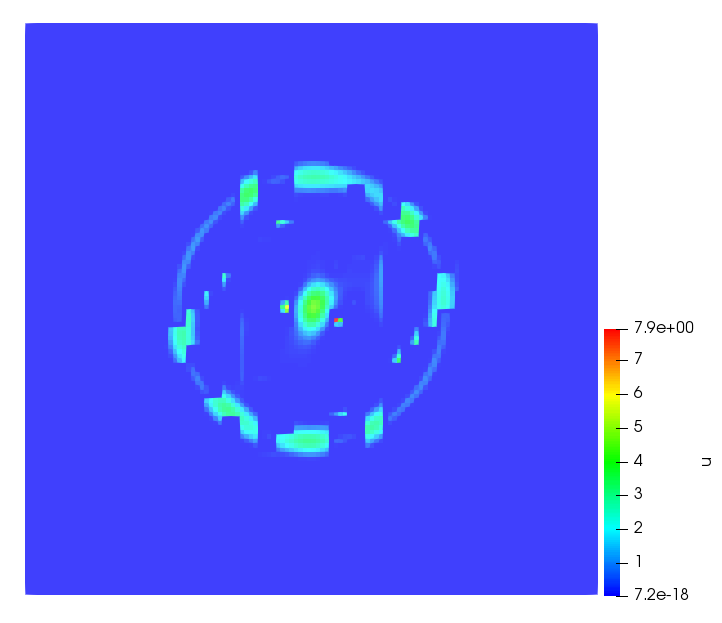}
		\quad \includegraphics[width=0.225\textwidth]{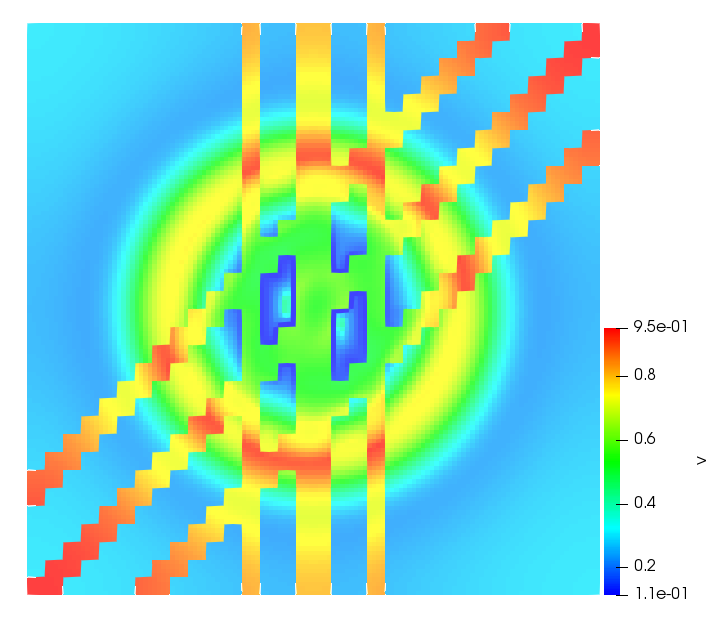}
		\quad \includegraphics[width=0.23\textwidth]{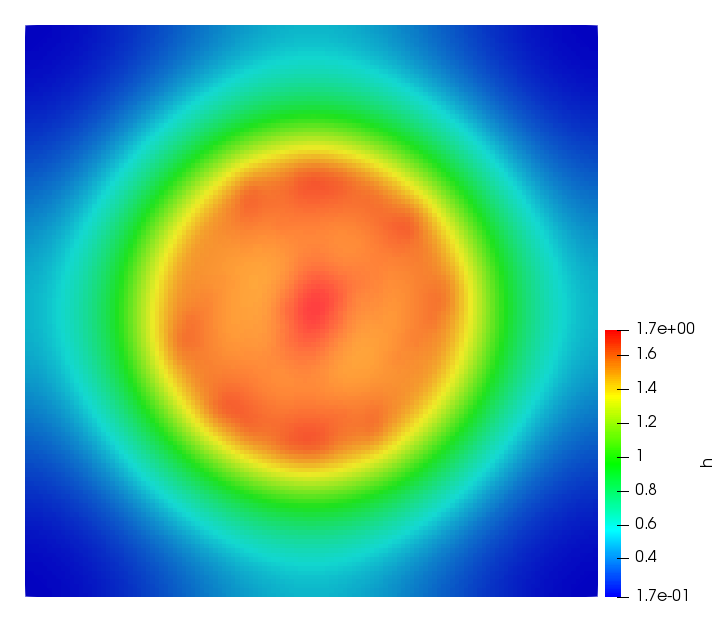}
		\quad \includegraphics[width=0.23\textwidth]{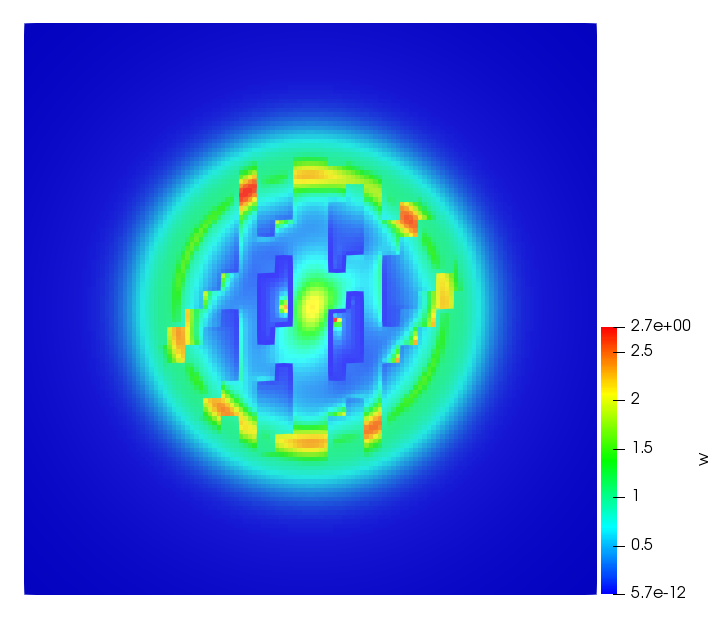}\\
		\caption{Evolution of tumor (first column), tissue (2nd column), MDEs (3rd column), and CAFs (last column) for model \eqref{eq:cafs} with $\mu =0$. Succesive times from top to bottom, top row: initial conditions.}\label{fig:cafs}
\end{figure}
\begin{figure}[h]
		\includegraphics[width=0.25\textwidth]{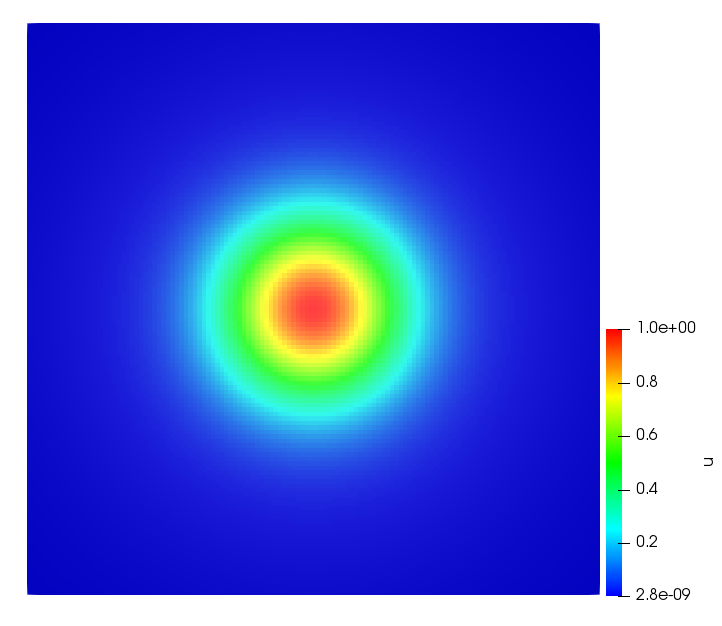}
		\quad \includegraphics[width=0.25\textwidth]{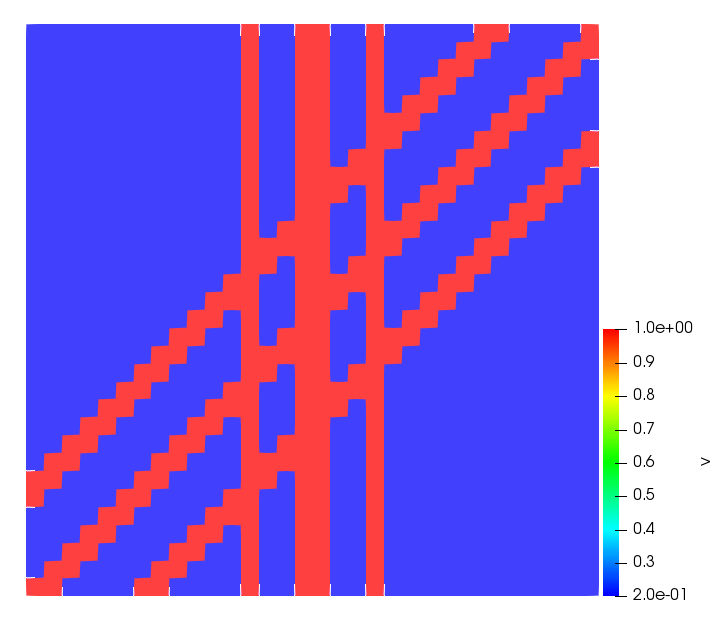}
		\quad \includegraphics[width=0.25\textwidth]{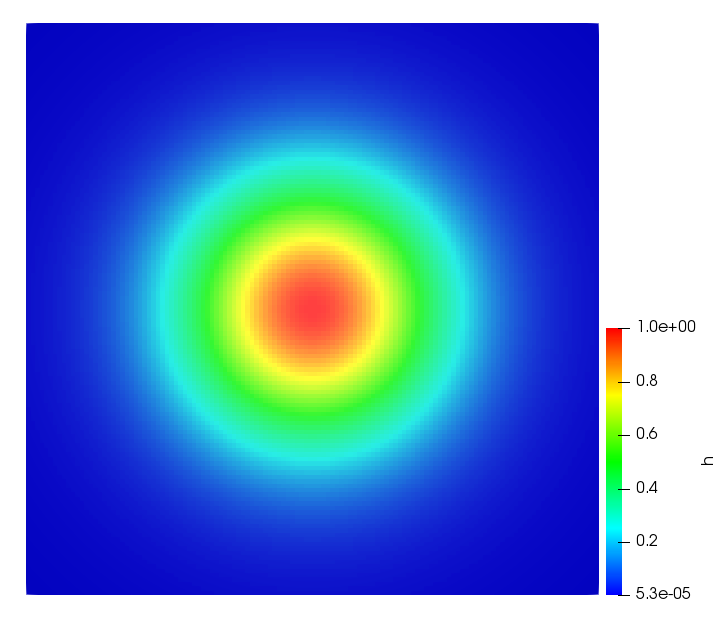}\\[1ex]
		\includegraphics[width=0.25\textwidth]{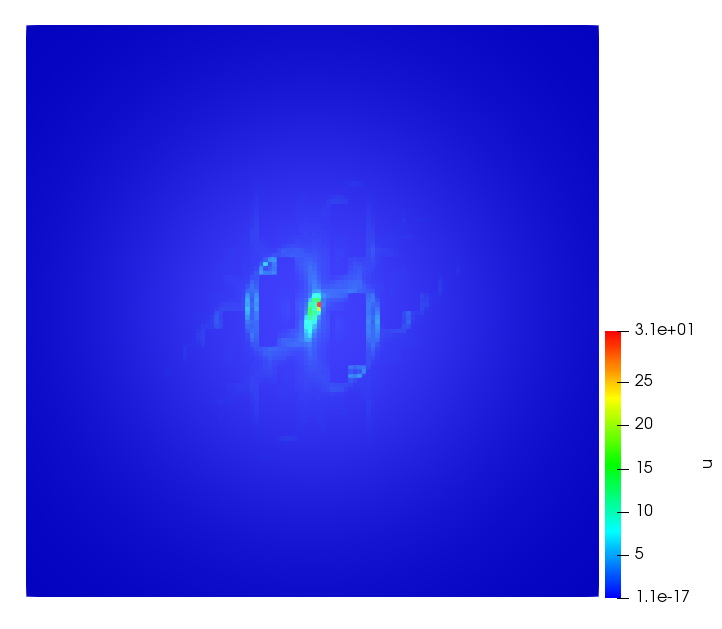}
		\quad \includegraphics[width=0.25\textwidth]{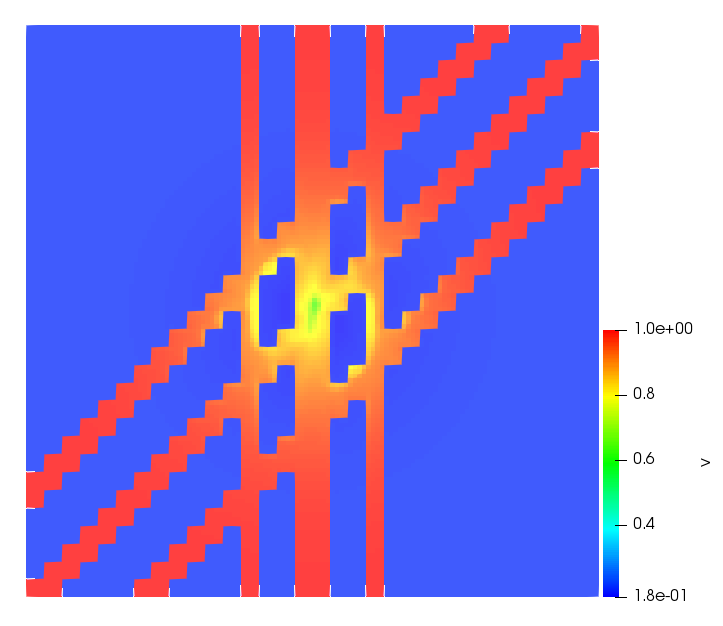}
		\quad \includegraphics[width=0.25\textwidth]{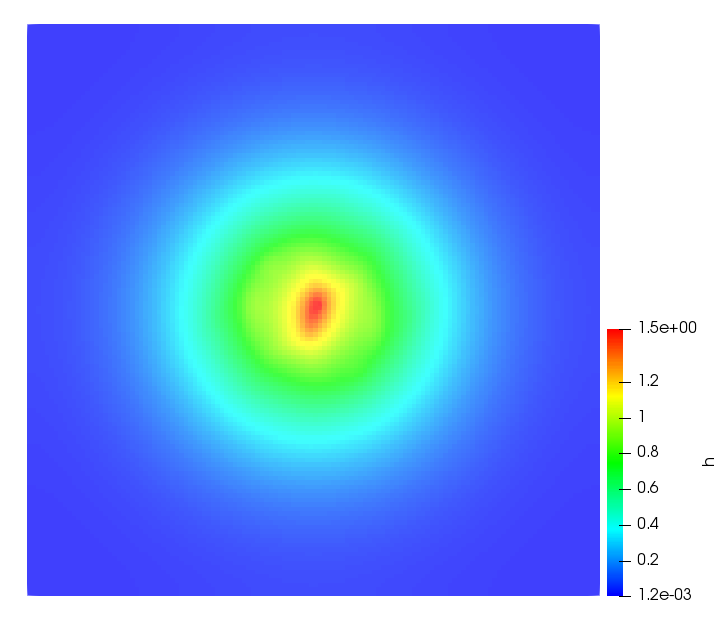}\\[1ex]
		\includegraphics[width=0.25\textwidth]{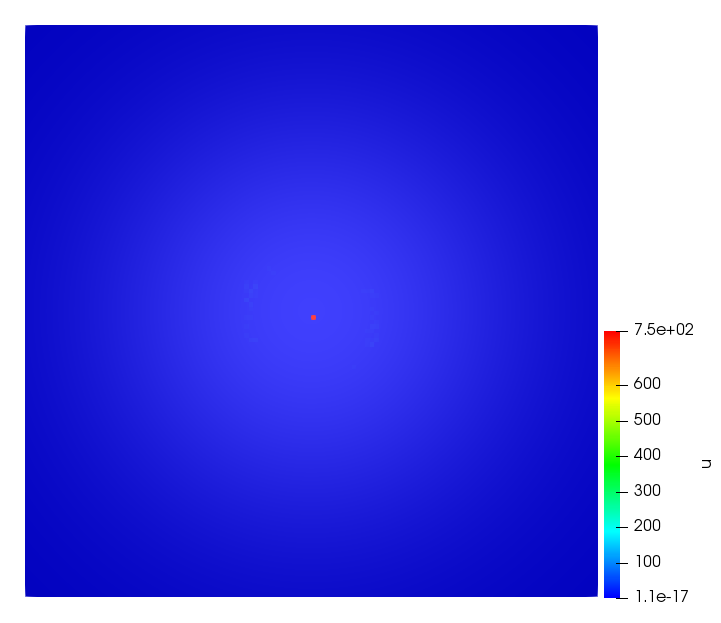}
		\quad \includegraphics[width=0.25\textwidth]{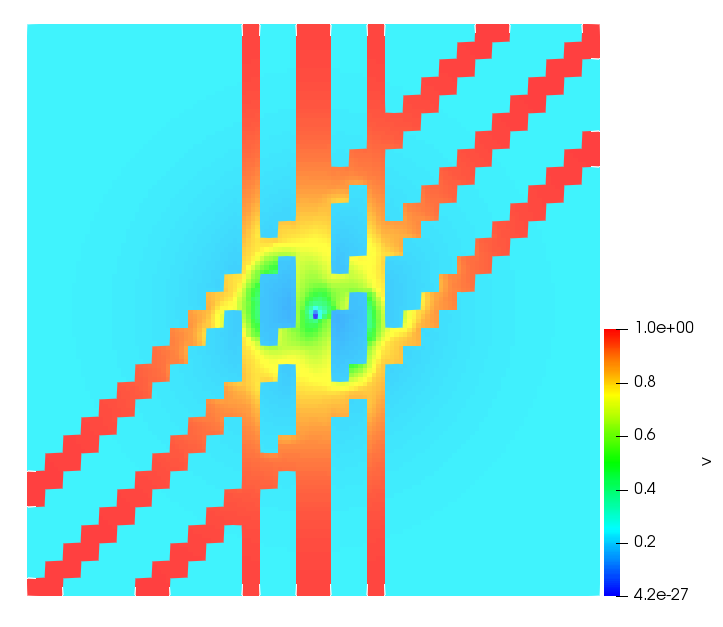}
		\quad \includegraphics[width=0.25\textwidth]{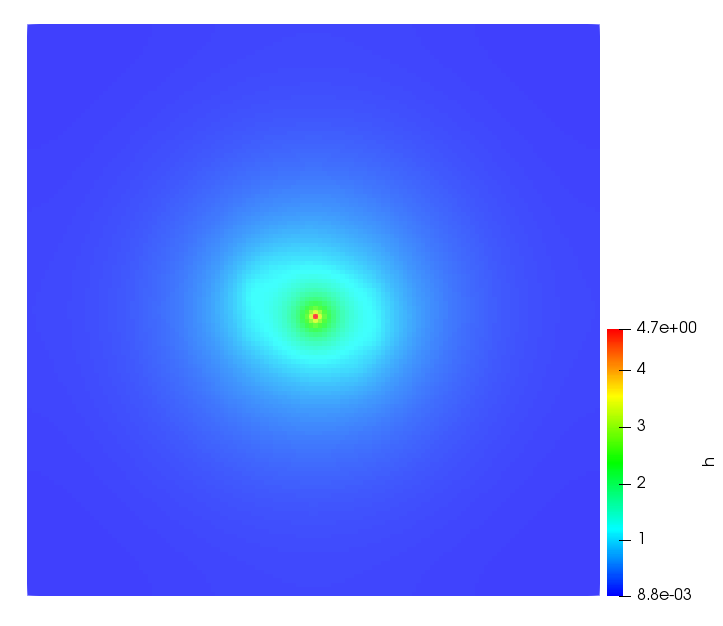}\\[1ex]
		\includegraphics[width=0.25\textwidth]{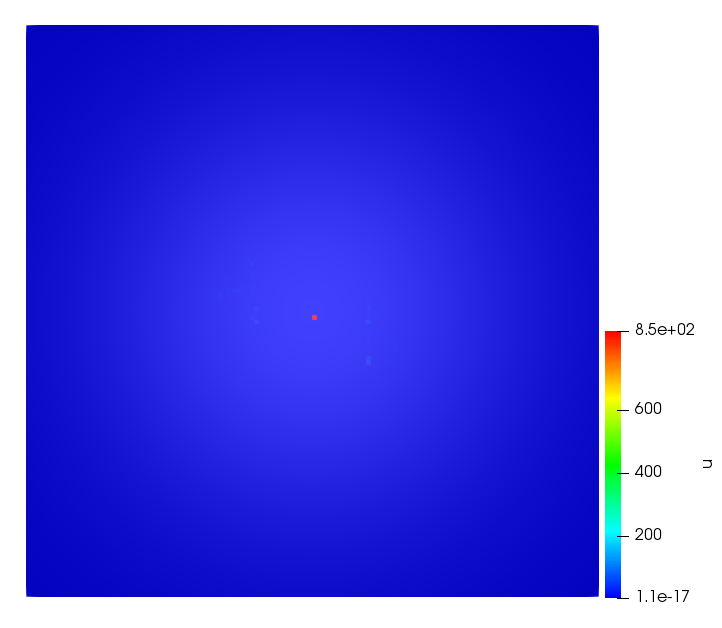}
		\quad \includegraphics[width=0.25\textwidth]{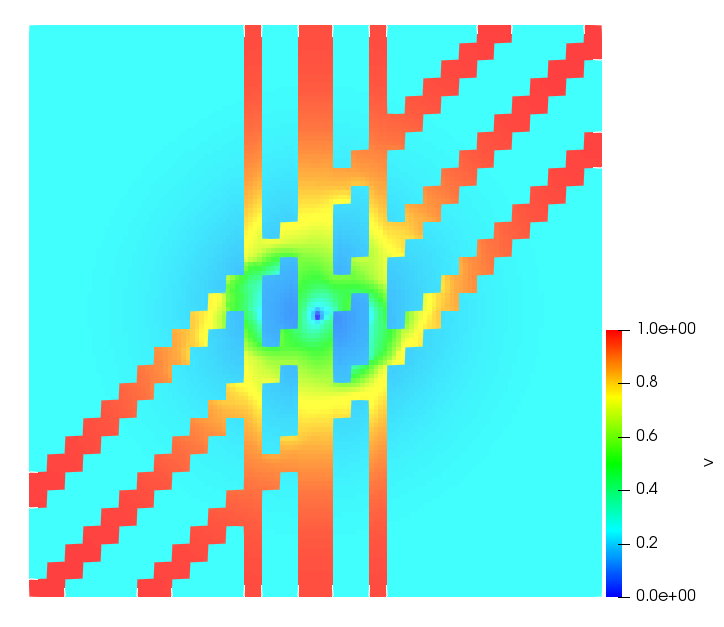}
		\quad \includegraphics[width=0.25\textwidth]{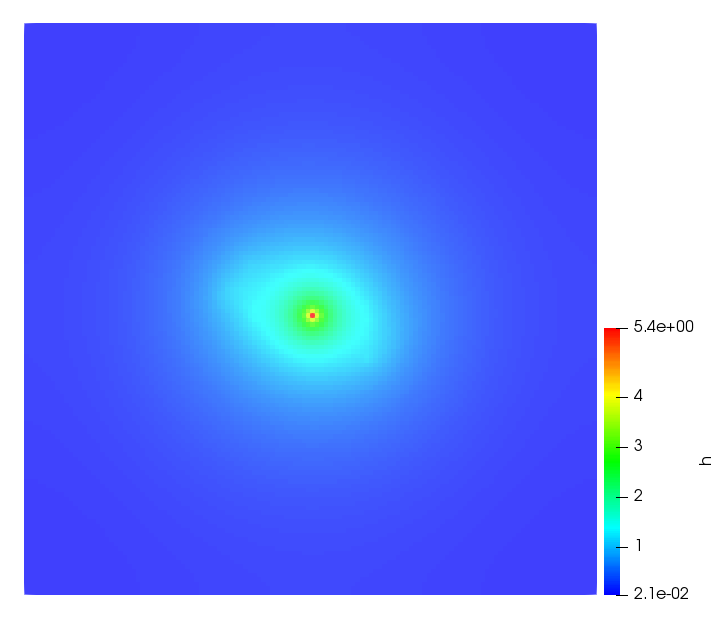}\\[1ex]
		\includegraphics[width=0.25\textwidth]{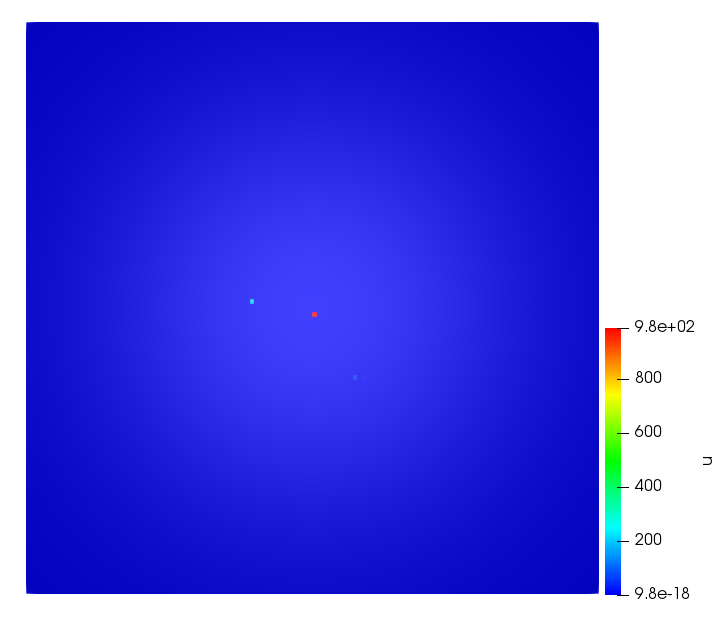}
		\quad \includegraphics[width=0.25\textwidth]{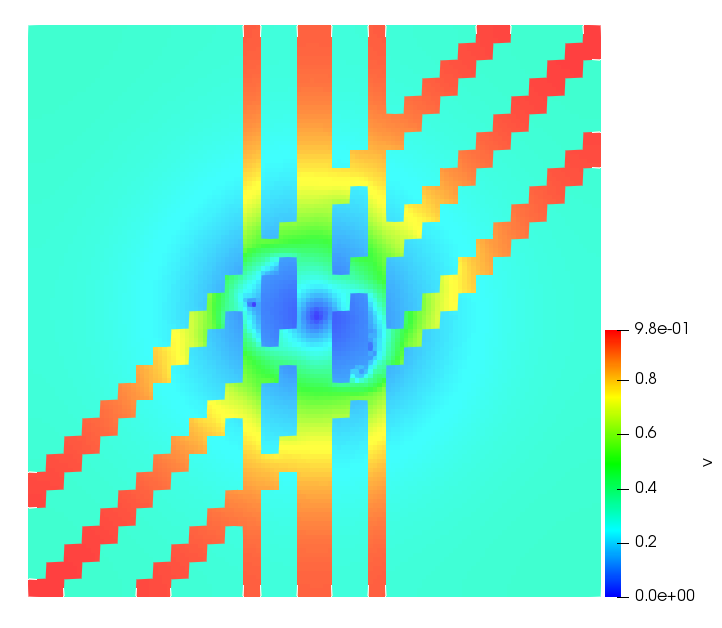}
		\quad \includegraphics[width=0.25\textwidth]{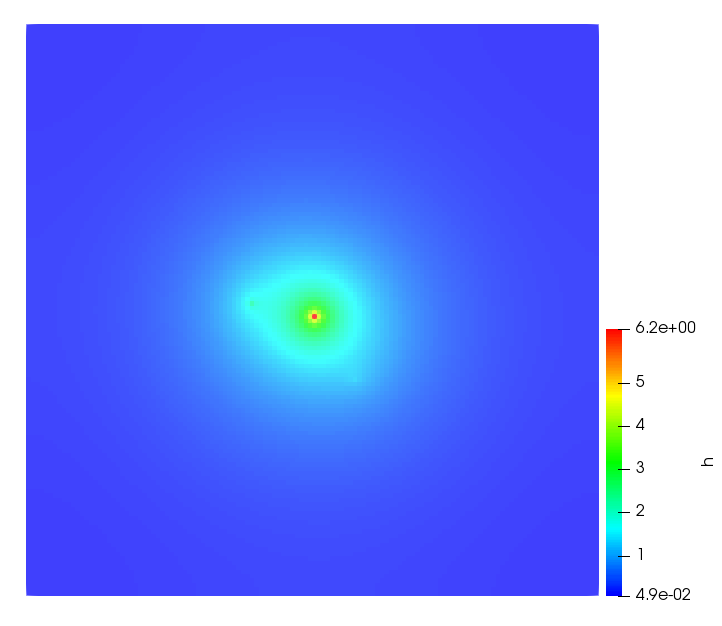}\\
		\caption{Evolution of tumor (first column), tissue (2nd column), and MDEs (last column) for model \eqref{eq:no-cafs} with $\mu =0$. Successive times from top to bottom, top row: initial conditions.}\label{fig:no-cafs}
\end{figure}

To illustrate our theoretical results we also performed some numerical simulations of the version
\be{eq:cafs}
\left\{ \begin{array}{l}
	u_t=D_u\Delta u - \chi \nabla \cdot (u\nabla h) - \xi \nabla \cdot (u\nabla v)+\mu u(1-u-v-w), \\[1mm]
	h_t=D_h\Delta h - h  + \alpha w, \\[1mm]
	v_t=-hv + \eta v(1-u-v)+\beta \frac{w}{1+w},\\[1mm]
	w_t=\gamma u,
	\end{array}\right .
\ee
of the indirect
signal production model \eqref{02} and -- for comparison purposes -- we also simulated some solutions of the system
\be{eq:no-cafs}
	\left\{ \begin{array}{l}
	u_t=D_u\Delta u - \chi \nabla \cdot (u\nabla h) - \xi \nabla \cdot (u\nabla v)+\mu u(1-u-v-w), \\[1mm]
	h_t=D_h\Delta h - h  + \alpha u, \\[1mm]
	v_t=-hv + \eta v(1-u-v),
	\end{array}\right .
\ee
in which the
signal (MDEs) is directly produced by the tumor cells. Since for both models there is no blow-up as long as $\mu >0$ 
(see \cite{zheng} for the result concerning \eqref{eq:no-cafs}), we only considered here the case with no tumor cell
proliferation, i.e. $\mu =0$. \abs
The simulations were performed using a discontinuous Galerkin FEM method. Thereby, the diffusion was discretized in space by
using a symmetric interior penalty Galerkin (SIPG) method (see \cite{riviere}), while for the drift term we did an upwind
discretization. The time was discretized with an IMEX procedure handling the diffusion implicitly and the reaction and taxis
terms explicitly. 
The computational domain was $\Omega =[0,1]^2$ and the initial conditions for tumor cells and MDEs were chosen in the form 
$u_0(x)=\exp(-\frac{|x|^2}{2\eps _u})$, $h_0(x)=\exp(-\frac{|x|^2}{2\eps _h})$, with $\eps _u=0.05$, $\eps_h=0.1$. For the
initial density of CAFs we considered a radially symetric form $w_0(r)=\exp (-\frac{(r-r_0)^2}{2\eps _w})$ with $\eps _w=0.01$
and $r_0=0.5$. Finally, the initial tissue density was characterized by $v_0(x)=\left \{\begin{array}{cc}
v_{\text{max}},&x\in \Omega _v\\                                                                                                                                                                                                                                                                                                                                                                                                                                                           
v_{\text{min}},&x\not\in \Omega _v                                                                                                                                                                                                                                                                                                                                                                                                                                                      \end{array}\right .,$
with $v_{\text{max}}=1$ and $v_{\text{min}}=0.2$, thus $\Omega _v\subset \Omega $ representing the stripes
shown in the second columns of Figures \ref{fig:cafs} and \ref{fig:no-cafs}. 
Together these initial conditions describe a heterogeneous tissue structure, in which a Gaussian-shaped tumor is embedded,
surrrounded by activated CAFs and featuring higher MDE concentration in the areas with many tumor cells. 
We considered for both models (where applicable) the following parameters: $\chi = 0.6$; $\xi = 0.5$; $D_u = 10^{-10}$; 
$D_h = 0.1$; $\eta = 10.6$; $\alpha = 5$; $\beta  = 1$; $\gamma = 1.0$. The results are shown in Figure \ref{fig:cafs} for
system \eqref{eq:cafs} and in Figure \ref{fig:no-cafs} for system \eqref{eq:no-cafs}. The first rows represent the initial
conditions, and the subsequent rows illustrate the solution behavior at several successive time steps. 

\newpage

The simulations of \eqref{eq:cafs} depicture an increase in tumor cell density, which is, however, limited. Same applies to the
MDE concentration and CAFs density, but with smaller rates. The tissue is correspondingly degraded, and the CAFs spread into
the region containing the main tumor mass, at the same time building up some tissue (for a sufficiently large $\beta $, as in
these simulations) at the sites where they are abundant enough. In contrast, model \eqref{eq:no-cafs} predicts localized tumor
cell aggregates of very high
density which are almost three orders of magnitude higher than the initial condition and keep growing in time, thus hinting on
blow-up of the solution. Likewise, in this latter setting the MDE concentration is directly produced by the tumor cells and
keeps growing as well, although to a much smaller rate than the cell density. The tissue degradation is much more localized 
and -- where it happens -- stronger. These results are conform with the theoretical findings in this 
and previous papers predicting blow-up of solutions when the signal was directly produced by the agents performing 
chemotaxis (\cite{BBTW}), while the solutions stay bounded in the case of indirect signal production.\abs
From a biological viewpoint the tumor cells use CAFs (which are originally 'harmless' stroma cells only becoming supporters of
tumor invasion upon activation) to produce matrix degrading factors. As mentioned above, the production of the latter seems to
be decisively controled by CAFs, hence is rather indirect, as the neoplastic cells first need to activate the CAFs, which then
enhance degradation of surrounding tissue and cell motility, including chemotaxis towards the gradient of proteolytic agents.
Our mathematical result actually tells that such mediation of invasion leads to avoidance of blow-up, unlike previous models
where the direct production of chemoattractant let the solution(s) become unbounded, a rather unrealistic biological scenario.
The result 
is in line with many in vivo and in vitro observations that cancer cells 'hijack' their environment in order to gain migratory,
survival, and proliferative advantages.\abs
We considered here that CAFs were non-diffusing, although they are indeed able to spread \cite{tao-caf}. Accounting for
diffusion of the chemoattractant producer $w$ does not pose, however, 
any further challenge to the analysis of our model \eqref{02}; in that case our setting belongs to the same mathematical class
as the one in \cite{chaplain-lolas}, to which it is also biologically 
related: both describe the evolution of a tumor under chemotaxis and haptotaxis, the chemoattractant(s) -- of which MDEs are
considered in both models -- being supposed to diffuse. \abs
Further chemotaxis-haptotaxis models belonging to the class studied here can be considered, of which \eqref{01} is just one 
example. 
As mentioned above, a diffusing signal producer can be easily accomodated to this model class -- if the diffusion is linear.
Solution-dependent diffusions of either involved species need further investigation; in \cite{QMW} one such system extending
that from \cite{taowin_JEMS} to allow for nonlinear diffusion of the chemotactic species has been studied, however without
haptotaxis.

\subsection*{Acknowledgement} 
The authors thank Aydar Uatay (TU Kaiserslautern) for performing the numerical simulations.
The second author acknowledges support by Deutsche
Forschungsgemeinschaft in the framework of the project  {\em
	Emergence of structures and advantages in cross-diffusion systems}
(No.~411007140, GZ: WI 3707/5-1).


%
%
%
%

%
%
%
%
%
%
%
%
%
%
%
%
\end{document}